\documentclass[a4paper,10pt]{amsart}

\makeatletter
\@namedef{subjclassname@2020}{\textup{2020} Mathematics Subject Classification}
\makeatother

\usepackage[latin1]{inputenc}
\usepackage{amssymb,amsmath}

\usepackage[all]{xy}

\usepackage[pagebackref]{hyperref}
     \hypersetup{
         colorlinks=true,
         linkcolor=cyan,
         citecolor=cyan
                     }
                     
\usepackage{lineno} 
 
\usepackage{mathtools}
\usepackage{tikz} 

\usepackage{color}

\newcommand{\cyan}[1]{\textcolor{cyan}{#1}}

\newtheorem{theorem}{Theorem}[section]

\newtheorem{proposition}[theorem]{Proposition}

\theoremstyle{definition}
\newtheorem{definition}[theorem]{Definition}

\newtheorem{remark}[theorem]{Remark}

\numberwithin{equation}{section}

\newcommand{\gras}[1]{{\mathbb #1}} 
\newcommand{\B}{\gras{B}}
\newcommand{\K}{\gras{K}}

\newcommand{\Z}{\gras{Z}} 
 
\newcommand{\R}{\gras{R}} 
\newcommand{\C}{\gras{C}}

\newcommand{\mult}{\mathfrak{m}}

\title{The Euclidean algorithm, lotuses and singularities}
\author[P. Popescu-Pampu]{Patrick Popescu-Pampu} 
\thanks{{\em Acknowledgments.} I acknowledge the support of the CDP C2EMPI, 
       together with the French State under the France-2030 programme, the University of Lille, 
       the Initiative of Excellence of the University of Lille, the European Metropolis of Lille for 
       their funding and support of the R-CDP-24-004-C2EMPI  project.  This research was funded 
       in part by l'Agence Nationale de la Recherche (ANR), project 
       ANR-22-CE40-0014 SINTROP.  I am grateful to Bernard Teissier for having recommended me  
       Fowler's book when I was his PhD student, to Evelia Garc\'{\i}a Barroso and 
      Pedro Gonz\'alez P\'erez for our long and fruitful collaboration around lotuses of plane 
      curve singularities, and to Jessica Carter and Cecilia Neve Jim\'enez for our philosophical 
      and historical discussions. Moreover, I thank them all, as well as Karine Chemla, 
      Isaac Gonz\'alez Rodr\'{\i}guez 
      and Paul-Emmanuel Timotei, for their comments on previous versions of this text.}
\date{26 May 2026}


\begin{document}

\begin{abstract}
        The {\em anthyphairetic process} leads from a pair $(a,b)$ of coprime positive integers 
        to the pair $(1,1)$ by successive subtractions of the smaller number from the bigger one. 
        This process, which is a slow version of Euclid's algorithm applied to the pair $(a,b)$, 
        corresponds naturally to the process of successive blowups leading to the minimal 
        embedded resolution 
        of the plane curve defined by $y^a - x^b = 0$.  This blowup 
        process may be represented graphically by a special two-dimensional simplicial 
        complex called a {\em lotus}. This allows to localize the various numbers appearing 
        either during the anthyphairetic process or during the Euclidean algorithm 
        at precise positions inside the lotus. In this introductory article, I recall first 
        the construction of this lotus starting 
        from the sequence of quotients generated by the Euclidean algorithm.  
        I present then an alternative way of constructing it directly from the sequence of 
        pairs of coprime integers generated by the anthyphairetic process, using what I call 
        {\em anthyphairetic rectangles}. I conclude by explaining how to reconstruct  from a lotus 
        the corresponding sequence of pairs of coprime integers. 
        This is a simple illustration of the way lotuses may serve as 
        {\em computational architectures}. 
\end{abstract}

{\bf This paper will appear in {\em Reading and Sharing Texts in the History of Science: An International Tribute to Karine Chemla}. F. Bretelle-Establet, M. C. Bustamante, E. Haffner, A. Keller, E. Sammarchi and X. Wang editors. Collection {\em Trends in the History of Science}, Springer. It is not in final form, therefore remarks are welcome!}
\bigskip

\maketitle

{\em \hfill I dedicate this article to Karine Chemla, 
   who likes colors in mathematical drawings. \smallskip}

\tableofcontents

\medskip
\section{Introduction}  \label{sect:intro}
\medskip

{\em Anthyphairesis} is a term classically used in the research on ancient Greek theories of 
ratio and proportion, related to what is called the {\em Euclidean algorithm} 
(see \cite{Cl 04} and \cite{K 95}). Here we 
will use the term in a sense close to David Fowler's \cite[Chapter 2]{F 99} or \cite[Section 6]{F 79}.  Let me quote from the introduction of \cite[Chapter 2]{F 99}:  

    \begin{quote}
        {\em I shall be arguing that substantial portions of the {\em Elements}, Plato's mathematical references, and other fourth- and third-century testimonies can be interpreted very plausibly in anthyphairetic terms. The sequence of repetition numbers that arises from anthyphairesis represents the result of repeated subtractions, not divisions, and this interpretation is corroborated by the
Greek name, derived from {\em anti-hypo-hairesis}, `reciprocal sub-traction'.}     
     \end{quote}

Let us say that an {\em anthyphairetic process} 
starts from two quantities $(A,B)$ and replaces them iteratively by new pairs of quantities obtained by subtracting each time the smaller one from the greater one. If this {\em anthyphairetic process} leads to a pair $(D,D)$ of equal quantities -- or, equivalently, in modern terms which accept the existence of a zero quantity, if they lead to a pair $(D,0)$ in which one of the quantities is zero -- then $A$ and $B$ are called {\em commensurable}; otherwise, they are called {\em incommensurable}.  In the commensurable case, the quantity $D$ is the {\em greatest common divisor} of $A$ and $B$.  Then $(A,B)$ are proportional to a unique pair $(a,b)$ of coprime positive integers, related to $A$ and $B$ by the equalities $A = a D$ and $B = b D$. 

Nowadays, the previous process of computation of the greatest common divisor is usually not presented as a finite process of iterative {\em subtractions}, but of iterative {\em divisions} with remainders, starting from the division of $B$ by $A$ (assuming that $B > A$). One speaks then of the {\em Euclidean algorithm}.  If $q_1, q_2, \dots, q_m$ are the positive integers appearing as quotients of the divisions performed during the algorithm, then in the commensurable case 
one gets a finite continued fraction expression of the ratio $B/A = b/a$:
    \[ \frac{b}{a} = q_1+\cfrac{1}{q_2+\cfrac{1}{ \cdots +\cfrac{1}{q_m}}}.  \]

  There is a strong relation between the anthyphairetic process and the theory of singularities 
  of plane curves.  Let us explain its simplest manifestation. Associate to each pair $(a,b)$ 
  of coprime positive integers the {\em binomial curve} $C_{a,b}$ defined by the equation 
  $y^a - x^b =0$. If $a, b > 1$, then this curve has a unique singular point, which lies 
  at the origin of coordinates. By blowing up this point, one transforms $C_{a,b}$
  either into the curve $C_{a,b-a}$ if $a < b$, or into the curve $C_{a-b,b}$ if $a > b$ 
  (see Proposition \ref{prop:prefchart}). If one 
  keeps blowing up the singular point of each new curve, one reaches therefore the curve 
  $C_{1,1}$ through a sequence of curves determined by the pairs of positive integers appearing 
  during the anthyphairetic process applied to the pair $(a,b)$. After one more blowup,  
  one gets the so-called {\em minimal embedded resolution} of the starting curve $C_{a,b}$. 
  
  Therefore, the process of minimal embedded resolution of $C_{a,b}$ by blowups of points 
  mimicks the anthyphairetic process applied to $(a,b)$. This fact explains 
  the appearance of the Euclidean algorithm and of continued fraction expansions of positive 
  rational numbers  in the study of plane curve singularities. Note that a second kind of 
  continued fraction expansions, the {\em Hirzebruch-Jung expansions}, 
  appears in the study of surface singularities (see  \cite{PP 07}). 
  
  Max Noether showed in \cite{N 90} how Euclidean algorithms appear naturally when one tries to 
  analyse the structure of a curve at one of its singular points by a sequence of blowups. 
  His description was extended by Federigo 
  Enriques and Oscar Chisini in  \cite[Libro IV, Capitolo I]{EC 18} to the case of germs of 
  plane curves admitting various branches. They introduced 
  in this context the first graphical representation 
  of such a blowup process, using what came to be known as {\em Enriques diagrams} 
  (see \cite[Section 3.9]{C 00} or \cite[Section 3]{PP 11}), which are special types of 
  rooted trees. 
  
  Later on, a second type of graphical representation of such a process was introduced, 
  namely the {\em dual graph of the 
  final exceptional divisor} (see Definition \ref{def:dualgraph}). One may 
  consult  \cite[Section 8.4]{BK 86}, \cite[Appendix]{LMW 89}, 
  \cite[Sections 5.3--5.4]{DP 00} or \cite[Section 3.6]{W 04} 
  for a description of their structure for plane curve singularities and \cite{PP 22} 
  for a historical discussion of the use of dual graphs in the study of surface singularities.  
  Both the Enriques diagram and the dual graph of a blowup process are finite trees 
  endowed with supplementary structures: straight/curved edges and broken vertices for 
  the Enriques diagram, negative vertex weights for the dual graph. 
  In general, they are not isomorphic even as abstract trees. 
  
  In \cite{PP 11} I introduced a new graphical object allowing us to understand the relation between 
  the two trees: a special kind of simplicial complex of dimension two, which I called a {\em kite}. 
  Later, in \cite{GGP 20}  and \cite{GGP 26}, Evelia Rosa Garc\'{\i}a Barroso, 
  Pedro Daniel Gonz\'alez P\'erez and I studied many more  invariants of plane 
  curve singularities using a simpler kind of 
  simplicial complex of dimension two, which we called a {\em lotus}.  
  All lotuses may be obtained by gluing simpler lotuses, which are associated to positive 
  rational numbers. As explained in \cite[Section 1.5.2]{GGP 20} and \cite[Section 11]{GGP 26}, 
  the lotus $\Lambda(b/a)$ associated to the rational number $b/a$ 
  may be constructed by first expanding $b/a$ as a continued fraction 
  and by triangulating then a given triangle in a way specified by this continued fraction. 
  
  \medskip
  
\begin{figure}[ht!] 
\begin{center}
\begin{tikzpicture}[x=0.6cm,y=0.6cm]

 \begin{scope}[shift={(0,0)}]
     \draw [fill=yellow](0,0) -- (6,0)--(6,6)--(0,6)--cycle;
     \draw [thick, color=black](0,0) -- (6,0);
     \draw [thick, color=black](0,1) -- (6,1);
     \draw [thick, color=black](0,2) -- (6,2);
     \draw [thick, color=black](0,3) -- (6,3);
     \draw [thick, color=black](0,4) -- (6,4);
     \draw [thick, color=black](0,5) -- (6,5);
     \draw [thick, color=black](0,6) -- (6,6);

     \draw [thick, color=black](0,0) -- (0,6);
     \draw [thick, color=black](6,0) -- (6,6);

     \node[draw,circle,inner sep=1.5pt,fill=black, color=black] at (0,0) {};    
             \node [right, color=black] at (0,0.3) {$4$};
     \node[draw,circle,inner sep=1.5pt,fill=black, color=black] at (0,1) {};
              \node [right, color=black] at (0,1.3) {$4$};
     \node[draw,circle,inner sep=1.5pt,fill=black, color=black] at (0,2) {};
               \node [right, color=black] at (0,2.3) {$4$};
     \node[draw,circle,inner sep=1.5pt,fill=black, color=black] at (0,3) {};
                \node [right, color=black] at (0,3.3) {$1$};
     \node[draw,circle,inner sep=1.5pt,fill=black, color=black] at (0,4) {};
               \node [right, color=black] at (0,4.3) {$1$};
     \node[draw,circle,inner sep=1.5pt,fill=black, color=black] at (0,5) {};
               \node [right, color=black] at (0,5.3) {$1$};
     \node[draw,circle,inner sep=1.5pt,fill=black, color=black] at (0,6) {};

     \node[draw,circle,inner sep=1.5pt,fill=black, color=black] at (6,0) {};
              \node [left, color=black] at (6,0.3) {$11$};
     \node[draw,circle,inner sep=1.5pt,fill=black, color=black] at (6,1) {};
               \node [left, color=black] at (6,1.3) {$7$};
     \node[draw,circle,inner sep=1.5pt,fill=black, color=black] at (6,2) {};
                \node [left, color=black] at (6,2.3) {$3$};
     \node[draw,circle,inner sep=1.5pt,fill=black, color=black] at (6,3) {};
                \node [left, color=black] at (6,3.3) {$3$};
     \node[draw,circle,inner sep=1.5pt,fill=black, color=black] at (6,4) {};
                \node [left, color=black] at (6,4.3) {$2$};
     \node[draw,circle,inner sep=1.5pt,fill=black, color=black] at (6,5) {};
                 \node [left, color=black] at (6,5.3) {$1$};
     \node[draw,circle,inner sep=1.5pt,fill=black, color=black] at (6,6) {};
     
      \node [below, color=black] at (3, -0.5) {$R(4, 11)$};  
     
          \draw [->, very thick, color=magenta](8.5 ,4)   --  (10.5, 4);          
          \draw [->, very thick, color=magenta](10.5 ,2)   --  (8.5, 2);  
            
          \node [above, color=black] at (9.5, 4) {Section \ref{sect:contract}};            
           \node [below, color=black] at (9.5, 2) {Section \ref{sec:lotintorect}};
     
\end{scope}

  
   \begin{scope}[shift={(12,0)}]
       \draw [fill=yellow](0,0) -- (6,0)--(3,6)--cycle;
      \draw [thick, color=black](0,0) -- (6, 0);
     \draw [thick, color=black](0,0) -- (5.5, 1);
     \draw [thick, color=black](0,0) -- (5,2);
     \draw [thick, color=black](1.5, 3) -- (5,2);
     \draw [thick, color=black](1.5,  3) -- (4.33, 3.33);
     \draw [thick, color=black](1.5, 3) -- (3.66,  4.66);
  
     \draw [thick, color=black](0,0) -- (3,6);
     \draw [thick, color=black](6,0) -- (3,6);
    
      \node[draw,circle,inner sep=1.5pt,fill=black, color=black] at (6, 0) {};
     \node[draw,circle,inner sep=1.5pt,fill=violet, color=black] at (0,0) {};    
     \node[draw,circle,inner sep=1.5pt,fill=black, color=black] at (5.5, 1) {};
     \node[draw,circle,inner sep=1.5pt,fill=red, color=black] at (5,2) {};
     \node[draw,circle,inner sep=1.5pt,fill=blue, color=black] at (1.5, 3) {};
     \node[draw,circle,inner sep=1.5pt,fill=black, color=black] at (4.33, 3.33) {};
     \node[draw,circle,inner sep=1.5pt,fill=black, color=black] at (3.66, 4.66) {};
     \node[draw,circle,inner sep=1.5pt,fill=orange, color=black] at (3,6) {};
     
      \node [below, color=black] at (3, -0.5) {$\Lambda(11/4)$};  
     
      \node [below, color=black] at (6, -0.2) {$e_1$};
      \node [below, color=black] at (0, -0.2) {$e_2$};

  \end{scope}

 \end{tikzpicture}
\end{center}
\caption{The equivalence of the anthyphairetic rectangle $R(4,11)$ and of 
     the lotus $\Lambda(11/4)$.}  
   \label{fig:equiv}
    \end{figure}
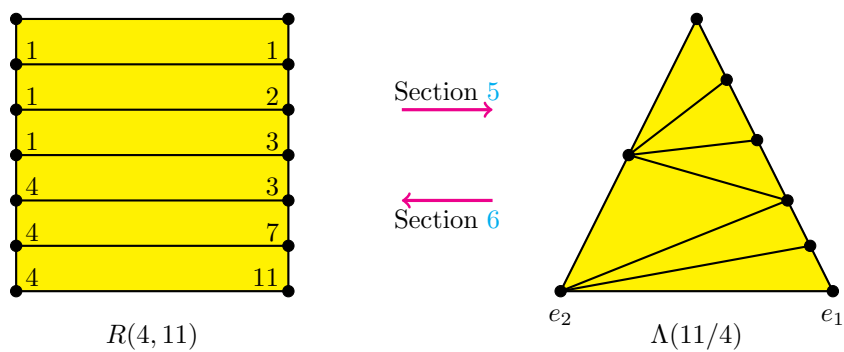

  In this paper I present an alternative way of constructing $\Lambda(b/a)$, starting from the 
  {\em anthyphairetic process} run on the coprime pair $(a,b)$ 
  (see Proposition \ref{prop:alternconstrlotus}).  
  Conversely, given a lotus, I explain how to get the pairs of positive integers 
  of the corresponding anthyphairetic process (see Proposition \ref{prop:revalg}). 
  Both constructions are based on a new way to represent graphically any 
  anthyphairetic process through an {\em anthyphairetic rectangle} and an associated 
  {\em anthyphairetic triangle} (see Definitions \ref{def:antrect} and \ref{def:triangpair}). 
   They are easier to understand on concrete examples than through general statements. 
  I chose to explain everything on the pair $(4, 11)$   (see Figure \ref{fig:equiv}), 
  as this case contains all the complexities of the general construction. 
  Therefore, it may be considered {\em paradigmatic},  
  in a sense similar to that used by Chemla in \cite{C 03} and Robadey in  \cite{R 04}.

  I tried to make this article accessible to people without prior knowledge of algebraic geometry or 
  singularity theory. In Section \ref{sect:blowuponce} I explain how a binomial 
  curve $C_{a,b}$ is transformed 
  by the blowup of the origin. In Section \ref{sect:embres} I describe the minimal embedded 
  resolution process of $C_{4, 11}$ by blowups of points, emphasizing how this 
  mimicks the anthyphairetic process which starts from the pair $(4, 11)$. 
  In Section \ref{sect:lotposrat} 
  I explain the notions of {\em universal lotus} and of {\em lotus of a positive rational number} 
  and I illustrate the method of  \cite[Section 1.5.2]{GGP 20} 
  and \cite[Section 11]{GGP 26} of construction of the lotus of a positive rational number 
  starting from its continued fraction expansion. In Section \ref{sect:contract} I present 
  an alternative construction of the lotus of a positive rational number $b/a$ 
  from the anthyphairetic process run on the pair $(a,b)$. Finally, in Section \ref{sec:lotintorect} 
  I explain how to reconstruct conversely from a lotus associated to a positive rational number 
  the sequence of pairs of coprime positive integers generated by the corresponding 
  anthyphairetic process. This provides a simple illustration of the way lotuses may serve as 
  {\em computational architectures}, a use amply studied in \cite{GGP 26}.

\medskip
\section{The blowup of an irreducible binomial curve}  \label{sect:blowuponce} 
\medskip

In this section I explain the standard representation of the blowup morphism of the origin 
of an affine plane by two monomial coordinate changes and how the curves defined 
by the binomial equations $y^a - x^b = 0$ are transformed by this blowup. 
\medskip

Let us fix a field $\K$ of characteristic zero (one may think for instance about $\R$ or $\C$).  
We will work inside smooth algebraic surfaces $S_i$ 
obtained from the affine plane $\boxed{\K^2_{x,y}}$ with coordinates $(x,y)$ 
by a finite sequence of blowups of points with coordinates in $\K$ (see Definition \ref{def:blowup}).   
We will explain below that such surfaces $S_i$ may be constructed by gluing finitely many copies 
of the affine plane using particular types of maps. We will call those copies {\bf charts}, 
because they endow $S_i$ with a structure of {\em atlas}, in the sense of differential geometry. 

We will look at irreducible and reduced curves contained in such a surface $S_i$. 
Such curves are subsets of $S_i$ whose intersections with each chart $\K^2_{u_j,v_j}$ 
is the vanishing set $Z(f_j)$ of a non-zero irreducible polynomial $f_j \in \K[u_j, v_j]$. 
More generally, by a {\bf curve} contained in such a smooth surface, which is then called the {\bf ambient surface} of the curve, we will mean a non-zero effective divisor, that is, a finite non-empty 
formal sum of irreducible and reduced curves with non-negative integer coefficients. 

If $f$ is a {\em rational function} defined on the smooth surface $S_j$ 
    (that is, a function which is a quotient of polynomials in each chart), 
   we denote by $\boxed{Z(f)}$ the curve defined by it, which is in this case {\em a principal divisor}. 
   For instance, if $m, n \in \Z$, then on the affine plane $\K^2_{x,y}$ one has the relation:
       \[ Z(x^m y^n) = m Z(x) + nZ(y), \]
   because the vanishing and the polar loci of the rational monomial 
   $x^m y^n \in \K[x,y]$ is included in the union 
   of the coordinate axes $Z(x)$ and $Z(y)$, and that the corresponding vanishing or polar 
   orders of $x^m y^n$ are $m$ and $n$ respectively. Namely, if $m \geq 0$, then one interprets 
   it as the vanishing order of $x^m y^n$ along $Z(x)$ and if $m <0$, then one interprets it as the 
   polar order of $x^m y^n$ along $Z(x)$. 
   
   Once one blows up 
   at least one point of $\K^2_{x,y}$, not all curves on the resulting surface are principal divisors. 
   For instance, the exceptional divisor of a blowup is not principal 
   (see point (\ref{it:oncebusi})  of Remark \ref{rem:selfint}).

\medskip
    The following computational description of {\em blowups} is crucial in order to establish 
 the parallel between the anthyphairetic processes and the resolution of binomial curves 
 by blowups:

\begin{definition}   \label{def:blowup}
    The surface $\boxed{\B_0 \K^2}$ obtained by blowing up the origin of $\K^2$ and  
    the associated {\bf blowup morphism} $\boxed{\pi} : \B_0 \K^2 \to \K^2$ 
    are obtained by gluing two copies of the affine plane $\K^2$, 
    with coordinates $(\bar{x}, \bar{y})$ and $(\tilde{x}, \tilde{y})$,  in such a way 
    that the morphism $\pi$ is described by the following changes of variables in those two charts:
   \begin{equation}     \label{eqn:blowupC2}
           \left\{ \begin{array}{l}      
                        x = \bar{x} \\
                        y = \bar{x} \bar{y}
                  \end{array} \right.   , \ \ \ 
           \left\{ \begin{array}{l}
                        x = \tilde{x} \tilde{y} \\
                        y = \tilde{y}
                  \end{array}  \right.  .
        \end{equation}
   \end{definition}     
        
   The gluing of the two affine planes $\K^2_{\bar{x}, \bar{y}}$ and $\K^2_{\tilde{x}, \tilde{y}}$ is 
   therefore described by the equations $\bar{x}  = \tilde{x} \tilde{y}$ and $\bar{x} \bar{y} = \tilde{y}$.  
   Those equations allow to express $(\tilde{x}, \tilde{y})$ as a rational function 
   of $(\bar{x}, \bar{y})$ and conversely. The corresponding rational maps 
      \[ (\bar{x}, \bar{y}) \mapsto (\bar{y}^{-1}, \bar{x} \bar{y}) \mbox{ and }  
          (\tilde{x}, \tilde{y}) \mapsto (\tilde{x} \tilde{y},  \tilde{x}^{-1}) \] 
    restrict to isomorphisms of their {\em determinacy loci} (the subsets of the corresponding 
    affine planes on which the two maps take well-defined values in $\K$), 
    which are $\K^2_{\bar{x}, \bar{y}} \setminus Z(\bar{y})$ 
   and $\K^2_{\tilde{x}, \tilde{y}} \setminus Z(\tilde{x})$.  They are mutually inverse {\bf birational maps}
       \[ \K^2_{\bar{x}, \bar{y}} \dashrightarrow   \K^2_{\tilde{x}, \tilde{y}},   \   
             \K^2_{\tilde{x}, \tilde{y}}  \dashrightarrow  \K^2_{\bar{x}, \bar{y}}, \]  
   that is, maps defined by rational fractions which establish an isomorphism between the 
   corresponding fields $\K(\bar{x}, \bar{y})$ and $\K(\tilde{x}, \tilde{y})$ of rational fractions. 
   Similarly, the blowup morphism $\pi$ is also birational. For more details about 
    blowups of points on smooth surfaces, see \cite[Section 8.4]{BK 86}, 
    \cite[Section 3.1]{C 00}, \cite[Section 3.2]{W 04} or \cite[Section 1.2.4]{GGP 20}. 
   
   In order to explain the relation between successive blowups and the anthyphairetic processes, 
   we need to emphasize the restrictions of $\pi$ to both charts $\K^2_{\bar{x}, \bar{y}}$ 
   and $\K^2_{\tilde{x}, \tilde{y}}$  of the surface $\B_0 \K^2$:
        
   \begin{definition} \label{def:affblowup}
         The morphism $\K^2_{\bar{x}, \bar{y}} \to \K^2_{x,y}$ described by the left-hand 
          system of (\ref{eqn:blowupC2}) is 
          called the {\bf $x$-blowup morphism} and the chart $\K^2_{\bar{x}, \bar{y}}$ is called 
          the {\bf $x$-chart of the blowup}. One defines analogously about the 
          {\bf $y$-blowup morphism} and the {\bf $y$-chart of the blowup}. 
   \end{definition}
   
   We will use the blowup morphism $\pi$ in order to {\em transform} curves $C$ passing 
   through the origin of $\K^2_{x,y}$. More generally, we will transform such curves by 
   compositions of blow ups. Such compositions are birational morphisms 
   $\varphi : S \to \K^2_{x,y}$ which restrict to isomorphisms between 
   $S \setminus \varphi^{-1}(0)$ and $\K^2_{x,y} \setminus \{0\}$. 
   There are two notions of {\em transform a curve $C \hookrightarrow \K^2_{x,y}$ by $\varphi$}: 
   
   \begin{definition}  \label{def:transforms}
       Let $\varphi : S \to \K^2_{x,y}$ be a birational morphism from a smooth surface $S$ 
       to $\K^2_{x,y}$, which restricts to an isomorphism between $S \setminus \varphi^{-1}(0)$ 
       and $\K^2_{x,y} \setminus \{0\}$. The {\bf exceptional divisor} of $\varphi$ 
       is the preimage $\varphi^{-1}(0)$ of the origin. Let $C$ be a curve lying on $\K^2_{x,y}$. 
       The {\bf total transform $\boxed{\varphi^*(C)}$ of $C$ by $\varphi$} is the divisor  
       $Z(\varphi^{*} (f))$ on $S$ defined by the {\bf pullback} $\boxed{\varphi^{*} (f)} := f \circ \varphi$  
       of a defining polynomial $f \in \K[x,y]$ of the curve $C$. The 
       {\bf strict transform $\boxed{\varphi^{-1}_s(C)}$ of $C$ by $\varphi$} 
       is the part of the total transform which consists of curves not contained in the exceptional divisor. 
   \end{definition}
   
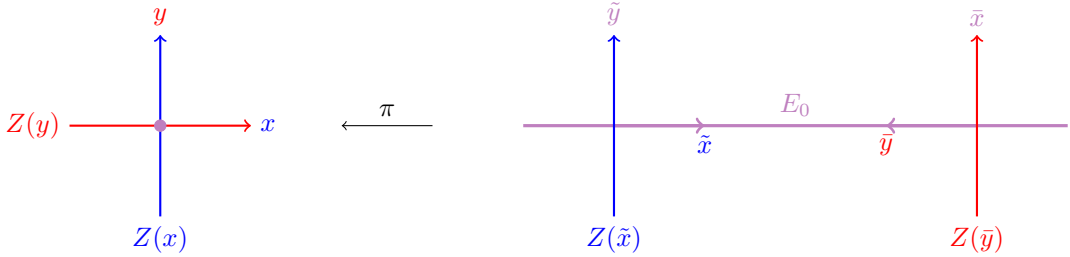
\begin{figure}[h!]
    \begin{center}
\begin{tikzpicture}[scale=1.2]

\begin{scope}[shift={(1,0)}]
        \draw [->, color=red, thick] (-1,0)--(1,0);
        \draw [->, color=blue, thick] (0,-1)--(0,1);

\node[draw,circle, inner sep=1.5pt,color=violet!50, fill=violet!50] at (0,0){};
\node [right, color=blue] at (1, 0) {$x$};
\node [below, color=blue] at (0, -1) {$Z(x)$};
\node [above, color=red] at (0, 1) {$y$};
\node [left, color=red] at (-1, 0) {$Z(y)$};
\end{scope}


\begin{scope}[shift={(8,0)}]
     \draw [-, color=violet!50, very thick](-3,0) -- (3,0);
           \draw [->, color=violet!50, very thick](-2,0) -- (-1,0);
           \draw [->, color=violet!50, very thick](2,0) -- (1,0);
     \draw [->, color=red, thick] (2,-1)--(2,1);
          \node [above, color=violet!50] at (-2, 1) {$\tilde{y}$};
          \node [below, color=red] at (2, -1) {$Z(\bar{y})$};
      \draw [->, color=blue, thick] (-2,-1)--(-2,1);
         \node [above, color=violet!50] at (2, 1) {$\bar{x}$};
         \node [below, color=blue] at (-2, -1) {$Z(\tilde{x})$};
         \node [below, color=blue] at (-1, 0) {$\tilde{x}$};
         \node [below, color=red] at (1, 0) {$\bar{y}$};
         \node [above, color=violet!50] at (0,0) {$E_0$};
\end{scope}


      \draw[<-](3,0)--(4,0);
           \node [above, color=black] at (3.5,0) {$\pi$};

\end{tikzpicture}
\end{center}
 \caption{The blowup morphism $\pi : \B_0 \K^2 = \K^2_{\tilde{x}, \tilde{y}} \cup \K^2_{\bar{x}, \bar{y}}  \to \K^2_{x,y}$ of $\K^2_{x,y}$ at the origin, its exceptional divisor $E_0$ and the strict transforms of the two axes of coordinates $Z(x)$ and $Z(y)$ of $\K^2_{x,y}$. A curve and its strict transform are 
 drawn with the same color.}
\label{fig:embrescusp}
   \end{figure}

   Denote by $\boxed{E_0}$ the {\bf exceptional divisor} of the blowup morphism 
   $\pi : \B_0 \K^2 \to \K^2$. 
   Formulae \eqref{eqn:blowupC2} show that $E_0$ is defined by the equations $\bar{x} =0$ 
   and $\tilde{y} =0$ respectively in the two charts. But if we look globally inside the union $\B_0 \K^2$ 
   of the two charts, we have the following total transforms of the axes of coordinates 
   $Z(x)$ and $Z(y)$ of $\K^2_{x,y}$ (see Figure \ref{fig:embrescusp}): 
        \begin{equation}     \label{eq:buaxes}
                \left\{ \begin{array}{l}      
                        \pi^*(Z(x)) =   Z(\tilde{x}) + E_0 ,   \\
                         \pi^*(Z(y)) = Z(\bar{y}) + E_0.
                  \end{array} \right.
        \end{equation}

   To each pair $(a,b) \in \Z^2_{>0}$ we associate a binomial and the corresponding 
   {\bf binomial curve}:   
     \[ \boxed{B_{a,b}} := y^a - x^b, \: \:   \boxed{C_{a,b}} := Z(y^a - x^b)  \hookrightarrow \K^2.  \]

    \begin{remark} 
       The curve $C_{a,b}$ is always {\em reduced}, because we assumed that the field $\K$ 
       is of characteristic zero (if it were of positive characteristic $p > 0$, then for instance 
       $y^p - x^p = (y-x)^p$ would define a non-reduced curve). Moreover, it is  {\em irreducible} 
       if and only if $a$ and $b$ are coprime (if $a = d a'$ and $b = d b'$ with $d > 1$, 
       then $B_{a',b'}$ would divide $B_{a,b}$ in the ring $\K[x,y]$, but not conversely). 
    \end{remark}

\medskip
{\bf In the sequel we will assume that $a$ and $b$ are coprime}.

Let us express the pullback $\pi^*(B_{a,b})$ of the binomial 
$B_{a,b}$ by the morphism $\pi$ in the two charts of $\B_0 \K^2$. 
{\em Assuming that $a < b$, we have:} 
      \begin{eqnarray*}
           B_{a,b}(\bar{x}, \bar{x} \bar{y}) =   (\bar{x} \bar{y})^a -  \bar{x}^b =  
                       \bar{x}^a ( \bar{y}^a -  \bar{x}^{b-a}),     \\
           B_{a,b}(\tilde{x} \tilde{y}, \tilde{y}) =   \tilde{y}^a - (\tilde{x} \tilde{y})^b =  
                      \tilde{y}^a(1 -   \tilde{x}^b \tilde{y}^{b-a}).
         \end{eqnarray*}  
One gets analogous formulae if $a > b$. 

We see that in both cases the total transform of $C_{a,b}$ decomposes as follows:
    \begin{equation}  \label{eq:decomptotstrict}
         \pi^*(C_{a,b}) =  \mult(C_{a,b}) E_0 + \pi^{-1}_s(C_{a,b}), 
    \end{equation}
 where $\boxed{ \mult(C_{a,b})}$ denotes the {\bf multiplicity} at the origin of the binomial curve
 $C_{a,b}$, that is, the smallest degree $\min\{a,b\}$ of the monomials appearing in its 
 defining polynomial $B_{a,b}$. Relation (\ref{eq:decomptotstrict}) holds in fact 
 for every curve $C$ contained in $\K^2_{x,y}$. 
 
 Note that {\em the strict transform $\pi^{-1}_s(C_{a,b})$ 
 is completely contained in exactly one of the two charts 
 of $\B_0 \K^2$}:  the unique chart whose origin belongs to $\pi^{-1}_s(C_{a,b})$. For this reason, 
 if one wants to understand this strict transform, then it is enough to work in that chart. That is:
    \begin{itemize}
       \item {\em If $a < b$}, it is enough to perform the $x$-blowup.  
            We get as strict transform of $C_{a,b}$: 
         \[
               \pi^{-1}_s(C_{a,b}) = Z( \bar{y}^a -  \bar{x}^{b-a}).
          \]
       \item {\em If $a > b$}, it is enough to perform the $y$-blowup.  
           We get as strict transform of $C_{a,b}$: 
         \[
               \pi^{-1}_s(C_{a,b}) = Z( \tilde{y}^{a-b} -  \tilde{x}^b).
          \]        
    \end{itemize} 

   As a consequence:

   \begin{proposition}   \label{prop:prefchart}  $\: $ 
      \begin{enumerate}
          \item Assume that the coprime positive integers $(a,b)$ verify the inequality $a < b$. 
             Then the strict transform of the binomial curve $C_{a,b}$ by the blowup of the origin 
             is isomorphic to the binomial curve $C_{a, b - a}$ and is contained in the $x$-chart. 
           \item Assume that the coprime positive integers $(a,b)$ verify the inequality $a > b$. 
             Then the strict transform of the binomial curve $C_{a,b}$ by the blowup of the origin 
             is isomorphic to the binomial curve $C_{a - b, b}$ and is contained in the $y$-chart.
        \end{enumerate}
   \end{proposition}

\medskip
\section{The minimal embedded resolution of an irreducible binomial curve}  \label{sect:embres}
\medskip

In this section I explain the relation between the anthyphairetic process applied to a 
pair $(a,b)$ of coprime positive integers and the process of minimal embedded resolution 
of the binomial curve $C_{a,b}$ by blowups of points. 
\medskip

Proposition \ref{prop:prefchart}  shows that at the level of pairs of exponents which determine an irreducible binomial curve, a blowup corresponds to an elementary subtraction in the {\em anthyphairetic process} analyzed by Fowler in \cite[Chapter 2]{F 99}. A finite process of blowups leads therefore to the binomial curve $C_{1,1}$. 

In order to illustrate this fact, let us consider the binomial curve
      \[  \boxed{C_{4,11} := Z(y^4 - x^{11})  \hookrightarrow \K^2_{x,y}}.  \] 
 We will blow up iteratively this curve and its strict transforms till reaching $C_{1,1}$. Instead of using the notations $(\bar{x}, \bar{y})$ and $(\tilde{x}, \tilde{y})$ of the changes of variables (\ref{eqn:blowupC2}), we will use indexed notations $(x_i, y_i)$. 
 
 Here $\boxed{\K^2_{x_i,y_i}}$ denotes the unique chart of the $i$-th blowup surface $\boxed{S_i}$ 
 which contains the whole strict transform of $C_{4,11}$. We indicate each change of coordinates, as well as the binomial defining polynomial of the strict transform in the new coordinate system: 
 \begin{equation}   \label{eq:seqbu}
   \begin{array}{l}
     y^4 - x^{11}  \rightsquigarrow 
     \left\{ \begin{array}{l}      
                        x = x_1 \\
                        y = x_1  y_1
                  \end{array} \right.    \rightsquigarrow 
       y_1^4 - x_1^7  \rightsquigarrow 
     \left\{ \begin{array}{l}      
                        x_1 = x_2 \\
                        y_1 = x_2  y_2
                  \end{array} \right.      \rightsquigarrow 
         y_2^4 - x_2^3  \rightsquigarrow 
     \left\{ \begin{array}{l}      
                        x_2 = x_3 y_3 \\
                        y_2 = y_3
                  \end{array} \right.      \rightsquigarrow     \\      $\ $ \\
       y_3^1 - x_3^3  \rightsquigarrow 
     \left\{ \begin{array}{l}      
                        x_3 = x_4 \\
                        y_3 = x_4  y_4
                  \end{array} \right.    \rightsquigarrow 
       y_4^1 - x_4^2  \rightsquigarrow 
     \left\{ \begin{array}{l}      
                        x_4 = x_5 \\
                        y_4 = x_5  y_5
                  \end{array} \right.  \rightsquigarrow 
       y_5^1 - x_5^1.
       \end{array}           
     \end{equation}
     
     The last binomial curve $Z(y_5^1 - x_5^1)$ is {\em smooth}. Moreover, the composition of 
     all the blowup morphisms is birational, therefore its restriction to $Z(y_5^1 - x_5^1)$ is a 
     birational morphism from the smooth curve $Z(y_5^1 - x_5^1)$ to the starting 
     curve $Z(y^4 - x^{11})$. 
     One says then that this restriction is a {\bf resolution of the singularities of $Z(y^4 - x^{11})$}. 
    Note that the two previous strict transforms $Z(y_3^1 - x_3^3)$ and $Z(y_4^1 - x_4^2)$ were 
    already smooth. Therefore, the restricted birational morphisms 
    $Z(y_3^1 - x_3^3) \to Z(y^4 - x^{11})$ and $Z(y_4^1 - x_4^2) \to Z(y^4 - x^{11})$ 
    were already resolutions of singularities of $Z(y^4 - x^{11})$. 
    
     Let us look now at those three strict transforms inside their ambient surfaces $ S_3, S_4, S_5$. 
     More precisely, we look at their positions relative to the exceptional divisors of the three 
     surfaces, which are the preimages of the origin of $\K^2_{x,y}$ by the corresponding 
     birational morphism $\pi_i : S_i \to \K^2_{x,y}$. 
     
     The curve $Z(y_3^1 - x_3^3)$ is 
     {\em tangent} to the component of the exceptional divisor described as $Z(y_3)$ in the chart 
     $\K^2_{x_3,y_3}$. Indeed, its tangent line at the origin is defined by the sum of the 
     terms of degree $1$ of the defining polynomial $y_3^1 - x_3^3$, 
     therefore it is the line $Z(y_3)$. Similarly, the curve $Z(y_4^1 - x_4^2)$ 
     is tangent to the component of the 
     exceptional divisor described as $Z(y_4)$ in the chart $\K^2_{x_4,y_4}$. By contrast, 
     the curve  $Z(y_5^1 - x_5^1)$ is {\em transversal} to both components $Z(x_5)$ and $Z(y_5)$ 
     of the exceptional divisor appearing in the chart $\K^2_{x_5,y_5}$.

      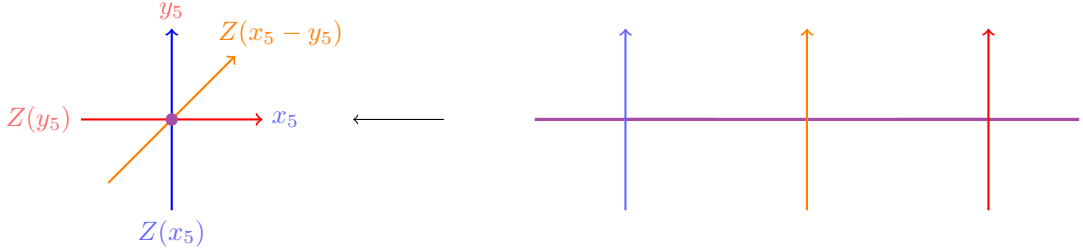
\begin{figure}[h!]
    \begin{center}
\begin{tikzpicture}[scale=1.2]

\begin{scope}[shift={(1,0)}]
        \draw [->, color=red, thick] (-1,0)--(1,0);
        \draw [->, color=blue, thick] (0,-1)--(0,1);
        \draw [->, color=orange, thick] (-0.7,-0.7)--(0.7,0.7);

\node[draw,circle, inner sep=1.5pt,color=violet!70, fill=violet!70] at (0,0){};
\node [right, color=blue!60] at (1, 0) {$x_5$};
\node [below, color=blue!60] at (0, -1) {$Z(x_5)$};
\node [above, color=red!60] at (0, 1) {$y_5$};
\node [left, color=red!60] at (-1, 0) {$Z(y_5)$};
\node [above, color=orange] at (1.2, 0.7) {$Z(x_5 - y_5)$};
\end{scope}


\begin{scope}[shift={(8,0)}]
     \draw [-, color=violet!70, very thick](-3,0) -- (3,0);
           \draw [-, color=violet!70, very thick](-2,0) -- (-1,0);
           \draw [-, color=violet!70, very thick](2,0) -- (1,0);
     \draw [->, color=red, thick] (2,-1)--(2,1);
      \draw [->, color=blue!60, thick] (-2,-1)--(-2,1);
       \draw [->, color=orange, thick] (0,-1)--(0,1);
\end{scope}


      \draw[<-](3,0)--(4,0);

\end{tikzpicture}
\end{center}
 \caption{The blowup morphism $\K^2_{x_5,y_5} \leftarrow \B_0 \K^2_{x_5,y_5}$ at the origin 
     of the chart $\K^2_{x_5,y_5} \hookrightarrow S_5$ 
     and the total transforms of the 
    curve $C_{4,11}$ both in $\K^2_{x_5,y_5}$ and in  $\B_0 \K^2_{x_5,y_5}$. 
    The strict transform of $C_{4,11}$ on the surface $S_5$ is the curve $Z(x_5 - y_5)$. 
    A curve on the left side $\K^2_{x_5,y_5}$ and its strict transform on the right side 
    are represented with the same color.}
\label{fig:threestrict}
   \end{figure}

     Composing the coordinate changes of (\ref{eq:seqbu}), we see that: 
         \[  \pi_5^*(B_{4, 11}) = \pi_5^*(y^4 - x^{11}) =  (x_5^8y_5^3)^4 - (x_5^3 y_5)^{11} 
            =  x_5^{32}  y_5^{12} - x_5^{33} y_5^{11}  = x_5^{32} y_5^{11} (y_5 - x_5). \] 
       This implies that in the chart $\K^2_{x_5,y_5}$ of $S_5$ the total transform of 
       the curve $B_{4, 11}$ is the divisor:
          \[ Z(\pi_5^*(B_{4, 11}))  =  32 \ Z(x_5) + 11 \ Z(y_5) + Z(y_5 - x_5). \]
       The strict transform of $B_{4, 11}$ is the component $Z(y_5 - x_5)$ of this divisor. 
     
     The three smooth irreducible curves $Z(x_5)$, $Z(y_5)$, $Z(y_5 - x_5)$ pass through 
     the origin of the chart (see Figure \ref{fig:threestrict}). 
     Therefore the total transform $Z(\pi_5^*(B_{4, 11}))$ does not look 
     like the union of coordinate axes of $\K^2$ in the neighborhood of this point: one says that    
     this total transform {\em is not a normal crossings divisor at that point}. Technically speaking, 
     this means that the underlying reduced curve of the divisor is not isomorphic under a formal 
     change of coordinates to the germ of the union of coordinate axes at one of their points.  
     By blowing up the origin of $\K^2_{x_5,y_5}$ we get a surface $S_6$ in which the total transform  
     $Z(\pi_6^*(B_{4, 11}))$ of $C_{4, 11}$ is a normal crossings divisor.  One says that 
     {\bf $\pi_6$ is an embedded resolution of the curve $C_{4, 11}$}. It is in fact the 
     {\bf minimal embedded resolution  of $C_{4, 11}$}, in the sense that any other embedded 
     resolution is obtained from $\pi_6$ by blowing up iteratively a finite set of points. 
     
     This fact illustrates the following general property, which goes back to Max Noether's \cite[Section I]{N 90} (see also \cite[Section 8.4, Theorems 12, 15]{BK 86}, \cite[Theorem 5.3.12]{DP 00}):
     
     \begin{proposition}   \label{prop:antires}
         Let $(a,b)$ be a pair of coprime positive integers. 
         The sequence of blowups leading to the penultimate stage of the minimal embedded resolution 
         of the binomial curve $C_{a,b}$ mimicks the sequence of subtractions of the anthyphairetic 
         algorithm applied to the pair $(a,b)$. 
     \end{proposition}
     
     Unlike the starting curve $C_{a,b}$ and the axes of coordinates $Z(x), Z(y)$, 
     the exceptional divisors created during the process of minimal embedded resolution 
     are {\em complete} in the sense of algebraic geometry (that is, they cannot be included 
     in larger irreducible curves in the way an affine line may be included in a projective line 
     by adding a point at infinity). Their 
     strict transforms are also complete. This implies that they have well-defined 
     {\em self-intersection numbers} in the ambient surfaces. These numbers may be computed 
     recursively, starting from the fact that at the moment of its creation, the self-intersection 
     of an exceptional divisor is $-1$ and that each time one blows up a smooth point of a 
     smooth complete curve, its self-intersection number drops by $1$ 
     (see Remark \ref{rem:selfint} or \cite[Lemma 8.1.6]{W 04}, \cite[Proposition 2.4]{GGP 26} 
     for more details).

      \begin{remark}   \label{rem:selfint}   
        $ \ $ 
        
    \begin{enumerate}   
       \item 
          The self-intersection number $A \cdot A \in \Z$ of a {\em complete}  
           irreducible curve $A$ may be defined by 
           taking a divisor $A'$ {\em linearly equivalent to $A$} (that is, such that 
           $A' = A + Z(f)$, where $f$ is a rational function on $S$), such that $A$ and $A'$ 
           have no common irreducible components, intersect transversally, 
           and by saying that $A\cdot A : = A \cdot A'$. 
           In turn, this number may be computed by imposing the bilinearity 
           of the intersection product and the convention that $B \cdot C$ is the number 
           of their intersection points, whenever $B$ and $C$ are irreducible curves which intersect 
           transversally. 
        \item 
           If the curve $A$ is not complete, then its self-intersection number 
           is not well-defined. For instance, consider $A := Z(y)$ in $\K^2_{x,y}$. By taking 
           successively $Z(y-1)$ and $Z(y-x)$ as divisor $A'$ in the construction above, 
           we get $0$ and $1$ as value of $A \cdot A'$.  
        \item  \label{it:oncebusi}
           A consequence of the definition above is that {\em the intersection number 
           of a complete curve with a principal divisor is $0$}. By applying this fact on the 
           blowup surface $S:= \B_0 \K^2$ of Definition \ref{def:blowup} to the exceptional 
          divisor $E_0$ (which is complete) and to the principal divisor $Z(\pi^* y)= Z(\overline{y}) + E_0$ 
          (see relations \ref{eq:buaxes}), we get $E_0 \cdot E_0 = -1$. As a consequence, 
          $E_0$ is not a principal divisor on $S$. 
        \item \label{it:dropsi}
           One may similarly prove that if one blows up a smooth point of a complete 
          irreducible curve $A$, then the self-intersection number of its strict transform $\tilde{A}$ 
          drops by $1$, that is, $\tilde{A} \cdot \tilde{A} = A \cdot A - 1$. Of course, 
          the ambient surfaces of $A$ and $\tilde{A}$ are distinct. 
    \end{enumerate}
   
 \end{remark}

     It is important to understand in terms of the pair $(a,b)$ 
     the structure of the total transform of the binomial curve $C_{a,b}$ by  
     its morphism of minimal embedded resolution. For instance, one may ask:  
     
        \medskip
      \begin{enumerate}
           \item[{\bf Q1)}] {\em What is the pattern of intersections of the various 
               irreducible components of the 
              exceptional divisor according to their order of appearance?} 
              
              \medskip
          \item[{\bf Q2)}] {\em What are their self-intersection numbers?} 
      \end{enumerate}
         \medskip
      
     It turns out that it is possible to understand {\em visually} the answers to those questions 
     and of similar 
     ones using a combinatorial object associated to the pair $(a,b)$, called its {\em lotus}. 
     In the next section we explain 
     its definition and we answer the two questions above (see Proposition \ref{prop:twoanswers}).

 \medskip
\section{The lotus of a positive rational number}  \label{sect:lotposrat}
\medskip

In this section I explain the definition of the lotus $\Lambda(b/a)$ of a positive rational number and 
its relation with the minimal embedded resolution process of the binomial curve defined by the equation $y^a - x^b =0$. 
\medskip

\begin{figure}[ht!]
     \begin{center}
\begin{tikzpicture}[scale=0.6]

\draw [fill=pink](1,0) -- (0,1)--(1,1)--cycle;
\draw [fill=pink!40](1,0) -- (1,1)--(2,1)--cycle;
\draw [fill=pink!40](1,0) -- (2,1)--(3,1)--cycle;
\draw [fill=pink!40](1,0) -- (3,1)--(4,1)--cycle;
\draw [fill=pink!40](1,0) -- (4,1)--(5,1)--cycle;
\draw [fill=pink!40](1,0) -- (5,1)--(6,1)--cycle;
\draw [fill=pink!40](1,0) -- (6,1)--(7,1)--cycle;
\draw [fill=pink!40](1,0) -- (7,1)--(8,1)--cycle;
\draw [fill=pink!40](1,0) -- (8,1)--(9,1)--cycle;
\draw [fill=pink!40](1,0) -- (9,1)--(10,1)--cycle;
      \draw [fill=pink!40](1,0) -- (10,1)--(11,1)--cycle;

\draw [fill=pink!40](1,1) -- (2,1)--(3,2)--cycle;
\draw [fill=pink!40](1,1) -- (3,2)--(4,3)--cycle;
\draw [fill=pink!40](1,1) -- (4,3)--(5,4)--cycle;
\draw [fill=pink!40](1,1) -- (5,4)--(6,5)--cycle;
\draw [fill=pink!40](1,1) -- (6,5)--(7,6)--cycle;
\draw [fill=pink!40](1,1) -- (7,6)--(8,7)--cycle;
\draw [fill=pink!40](1,1) -- (8,7)--(9,8)--cycle;
\draw [fill=pink!40](1,1) -- (9,8)--(10,9)--cycle;
   \draw [fill=pink!40](1,1) -- (10,9)--(11,10)--cycle;
   
    \draw [fill=pink!40](11,2) -- (5,1)--(6,1)--cycle;
     \draw [fill=pink!40](11,3) -- (4,1)--(7,2)--cycle;
      \draw [fill=pink!40](11,8) -- (4,3)--(7,5)--cycle;
       \draw [fill=pink!40](11,9) -- (5,4)--(6,5)--cycle;
       
         \draw [fill=pink!40](2,11) -- (1,5)--(1,6)--cycle;
     \draw [fill=pink!40](3,11) -- (1,4)--(2,7)--cycle;
      \draw [fill=pink!40](8,11) -- (3,4)--(5,7)--cycle;
       \draw [fill=pink!40](9,11) -- (4,5)--(5,6)--cycle;

\draw [fill=pink!40](2,1) -- (3,2)--(5,3)--cycle;
\draw [fill=pink!40](2,1) -- (5,3)--(7,4)--cycle;
\draw [fill=pink!40](2,1) -- (7,4)--(9,5)--cycle;
   \draw [fill=pink!40](2,1) -- (9,5)--(11,6)--cycle;

\draw [fill=pink!40](2,1) -- (3,1)--(5,2)--cycle;
\draw [fill=pink!40](2,1) -- (5,2)--(7,3)--cycle;
\draw [fill=pink!40](2,1) -- (7,3)--(9,4)--cycle;
   \draw [fill=pink!40](2,1) -- (9,4)--(11,5)--cycle;

\draw [fill=pink!40](3,2) -- (5,3)--(8,5)--cycle;
   \draw [fill=pink!40](3,2) -- (8,5)--(11, 7)--cycle;

\draw [fill=pink!40](3,2) -- (4,3)--(7,5)--cycle;
\draw [fill=pink!40](3,2) -- (7,5)--(10,7)--cycle;

\draw [fill=pink!40](4,3) -- (5,4)--(9,7)--cycle;

\draw [fill=pink!40](3,1) -- (5,2)--(8,3)--cycle;
    \draw [fill=pink!40](3,1)--(8,3)--(11,4)--cycle;
\draw [fill=pink!40](3,1) -- (4,1)--(7,2)--cycle;
\draw [fill=pink!40](3,1) -- (7,2)--(10,3)--cycle;

\draw [fill=pink!40](4,1) -- (5,1)--(9,2)--cycle;

\draw [fill=pink!40](0,1) -- (1,1)--(1,2)--cycle;
\draw [fill=pink!40](0,1) -- (1,2)--(1,3)--cycle;
\draw [fill=pink!40](0,1) -- (1,3)--(1,4)--cycle;
\draw [fill=pink!40](0,1) -- (1,4)--(1,5)--cycle;
\draw [fill=pink!40](0,1) -- (1,5)--(1,6)--cycle;
\draw [fill=pink!40](0,1) -- (1,6)--(1,7)--cycle;
\draw [fill=pink!40](0,1) -- (1,7)--(1,8)--cycle;
\draw [fill=pink!40](0,1) -- (1,8)--(1,9)--cycle;
\draw [fill=pink!40](0,1) -- (1,9)--(1,10)--cycle;
   \draw [fill=pink!40](0,1) -- (1,10)--(1,11)--cycle;

\draw [fill=pink!40](1,1) -- (1,2)--(2,3)--cycle;
\draw [fill=pink!40](1,1) -- (2,3)--(3,4)--cycle;
\draw [fill=pink!40](1,1) -- (3,4)--(4,5)--cycle;
\draw [fill=pink!40](1,1) -- (4,5)--(5,6)--cycle;
\draw [fill=pink!40](1,1) -- (5,6)--(6,7)--cycle;
\draw [fill=pink!40](1,1) -- (6,7)--(7,8)--cycle;
\draw [fill=pink!40](1,1) -- (7,8)--(8,9)--cycle;
\draw [fill=pink!40](1,1) -- (8,9)--(9,10)--cycle;
   \draw [fill=pink!40](1,1) -- (9,10)--(10,11)--cycle;

\draw [fill=pink!40](1,2) -- (2,3)--(3,5)--cycle;
\draw [fill=pink!40](1,2) -- (3,5)--(4,7)--cycle;
\draw [fill=pink!40](1,2) -- (4,7)--(5,9)--cycle;
   \draw [fill=pink!40](1,2)--(5,9)--(6,11)--cycle;

\draw [fill=pink!40](1,2) -- (1,3)--(2,5)--cycle;
\draw [fill=pink!40](1,2) -- (2,5)--(3,7)--cycle;
\draw [fill=pink!40](1,2) -- (3,7)--(4,9)--cycle;
   \draw [fill=pink!40](1,2)--(4,9)--(5,11)--cycle;

\draw [fill=pink!40](2,3) -- (3,5)--(5,8)--cycle;
   \draw [fill=pink!40](2,3)--(5,8)--(7,11)--cycle;

\draw [fill=pink!40](2,3) -- (3,4)--(5,7)--cycle;
\draw [fill=pink!40](2,3) -- (5,7)--(7,10)--cycle;

\draw [fill=pink!40](3,4) -- (4,5)--(7,9)--cycle;

\draw [fill=pink!40](1,3) -- (2,5)--(3,8)--cycle;
     \draw [fill=pink!40](1,3)--(3,8)--(4,11)--cycle;
\draw [fill=pink!40](1,3) -- (1,4)--(2,7)--cycle;
\draw [fill=pink!40](1,3) -- (2,7)--(3,10)--cycle;

\draw [fill=pink!40](1,4) -- (1,5)--(2,9)--cycle;

\draw [dashed, gray] (0,0) grid (11,11);
\node [below] at (1,0) {$e_{1}$}; 
\node [left] at (0,1) {$e_{2}$}; 
\node [left] at (0,0) {$0$}; 

\draw[->][thick, color=black](-1,1) .. controls (-0.5,0.5) ..(0.7,0.7);  
\node [left] at (-1, 1) {$\delta(e_{1},e_{2})$}; 
   \end{tikzpicture}
\end{center}
  \caption{Partial view of the universal lotus $\Lambda(e_{1},e_{2})$ of Definition \ref{def:petalsunivlotus}: the union of its petals contained in the square $[0, 11] \times [0, 11]$.}  
  \label{fig:Unilotus} 
    \end{figure}
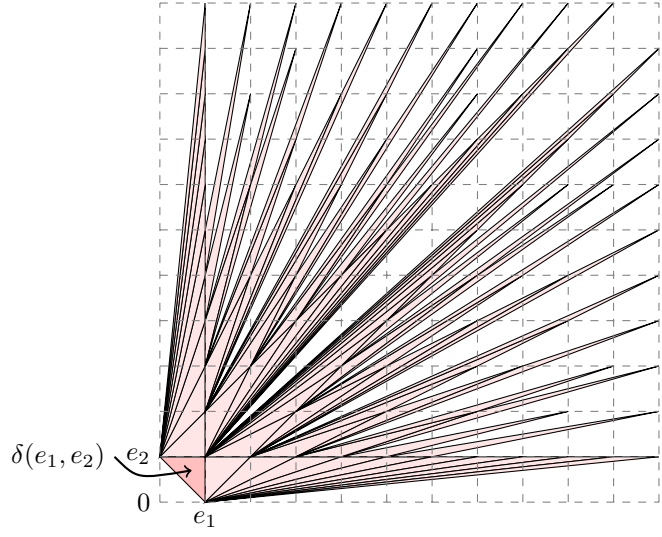

The following definition is a special case of  \cite[D\'efinition 5.2]{PP 11} and  
\cite[Definition 1.5.3]{GGP 20}:

\begin{definition}   \label{def:petalsunivlotus}
   Let $N$ be a free abelian group of rank two endowed with a basis $(e_1, e_2)$. The {\bf associated petal} of  $(e_1, e_2)$ is the triangle $\boxed{\delta(e_1, e_2)}$ with vertices $e_1, e_2, e_1 + e_2$ contained in the real plane $\boxed{N_{\R}} := N \otimes_{\Z} \R$. Its {\bf base} is the segment $[e_1, e_2]$ and its {\bf summit} is the vertex $e_1 + e_2$.  The {\bf universal lotus} $\boxed{\Lambda(e_1, e_2)}$ associated to $(e_1, e_2)$ (see Figure \ref{fig:Unilotus}) is the simplicial complex of pure dimension two contained in  $N_{\R}$ whose triangles are the petals associated to the bases of the lattice $N$ which have non-negative coordinates relative to the basis $(e_1, e_2)$. 
\end{definition}

One may associate a finite sublotus of the universal lotus to every positive rational number (see  \cite[D\'efinition 8.6]{PP 11} and \cite[Definition 1.5.4]{GGP 20}):

\begin{definition}   \label{def:lotusrat}
   Let $(a,b)$ be a pair of coprime positive integers.  The {\bf lotus} $\boxed{\Lambda(b/a)}$ 
   of $b/a$ (see Figure \ref{fig:Lotus11/4}) is the subcomplex of the universal lotus 
   $\Lambda(e_1, e_2)$ which is the set of petals whose interiors intersect 
   the segment $[0, a e_1 + b e_2]$ and of their edges and vertices. Its {\bf base} 
   is the segment $[e_1, e_2]$. 
\end{definition}

   As shown on the left side of Figure \ref{fig:Lotus11/4}, the lotus of a positive rational number 
  becomes rapidly very thin when $a$ and $b$ increase, which makes very difficult seeing 
  how its constituent petals are arranged inside it. That is why it is better to draw an isomorphic 
  but thicker simplicial complex (with two vertices decorated by the symbols $e_1$ and $e_2$, 
  in order to be able to extract from it the number $b/a$ without ambiguities), as on the right 
  side of Figure \ref{fig:Lotus11/4}.

 \begin{figure}[tb]
     \begin{center}
\begin{tikzpicture}[scale=0.6]

 \begin{scope}[shift={(0,0)}]

\draw [fill=yellow](1,0) -- (0,1)--(1,1)--cycle;
\draw [fill=pink!40](1,0) -- (1,1)--(2,1)--cycle;
\draw [fill=pink!40](1,0) -- (2,1)--(3,1)--cycle;
\draw [fill=pink!40](1,0) -- (3,1)--(4,1)--cycle;
\draw [fill=pink!40](1,0) -- (4,1)--(5,1)--cycle;
\draw [fill=pink!40](1,0) -- (5,1)--(6,1)--cycle;
\draw [fill=pink!40](1,0) -- (6,1)--(7,1)--cycle;
\draw [fill=pink!40](1,0) -- (7,1)--(8,1)--cycle;
\draw [fill=pink!40](1,0) -- (8,1)--(9,1)--cycle;
\draw [fill=pink!40](1,0) -- (9,1)--(10,1)--cycle;
      \draw [fill=pink!40](1,0) -- (10,1)--(11,1)--cycle;

\draw [fill=pink!40](1,1) -- (2,1)--(3,2)--cycle;
\draw [fill=pink!40](1,1) -- (3,2)--(4,3)--cycle;
\draw [fill=pink!40](1,1) -- (4,3)--(5,4)--cycle;
\draw [fill=pink!40](1,1) -- (5,4)--(6,5)--cycle;
\draw [fill=pink!40](1,1) -- (6,5)--(7,6)--cycle;
\draw [fill=pink!40](1,1) -- (7,6)--(8,7)--cycle;
\draw [fill=pink!40](1,1) -- (8,7)--(9,8)--cycle;
\draw [fill=pink!40](1,1) -- (9,8)--(10,9)--cycle;
   \draw [fill=pink!40](1,1) -- (10,9)--(11,10)--cycle;
   
    \draw [fill=pink!40](11,2) -- (5,1)--(6,1)--cycle;
     \draw [fill=pink!40](11,3) -- (4,1)--(7,2)--cycle;
      \draw [fill=pink!40](11,8) -- (4,3)--(7,5)--cycle;
       \draw [fill=pink!40](11,9) -- (5,4)--(6,5)--cycle;
       
         \draw [fill=pink!40](2,11) -- (1,5)--(1,6)--cycle;
     \draw [fill=pink!40](3,11) -- (1,4)--(2,7)--cycle;
      \draw [fill=pink!40](8,11) -- (3,4)--(5,7)--cycle;
       \draw [fill=pink!40](9,11) -- (4,5)--(5,6)--cycle;

\draw [fill=pink!40](2,1) -- (3,2)--(5,3)--cycle;
\draw [fill=pink!40](2,1) -- (5,3)--(7,4)--cycle;
\draw [fill=pink!40](2,1) -- (7,4)--(9,5)--cycle;
   \draw [fill=pink!40](2,1) -- (9,5)--(11,6)--cycle;

\draw [fill=pink!40](2,1) -- (3,1)--(5,2)--cycle;
\draw [fill=pink!40](2,1) -- (5,2)--(7,3)--cycle;
\draw [fill=pink!40](2,1) -- (7,3)--(9,4)--cycle;
   \draw [fill=pink!40](2,1) -- (9,4)--(11,5)--cycle;

\draw [fill=pink!40](3,2) -- (5,3)--(8,5)--cycle;
   \draw [fill=pink!40](3,2) -- (8,5)--(11, 7)--cycle;

\draw [fill=pink!40](3,2) -- (4,3)--(7,5)--cycle;
\draw [fill=pink!40](3,2) -- (7,5)--(10,7)--cycle;

\draw [fill=pink!40](4,3) -- (5,4)--(9,7)--cycle;

\draw [fill=pink!40](3,1) -- (5,2)--(8,3)--cycle;
    \draw [fill=pink!40](3,1)--(8,3)--(11,4)--cycle;
\draw [fill=pink!40](3,1) -- (4,1)--(7,2)--cycle;
\draw [fill=pink!40](3,1) -- (7,2)--(10,3)--cycle;

\draw [fill=pink!40](4,1) -- (5,1)--(9,2)--cycle;

\draw [fill=yellow](0,1) -- (1,1)--(1,2)--cycle;
\draw [fill=yellow](0,1) -- (1,2)--(1,3)--cycle;
\draw [fill=pink!40](0,1) -- (1,3)--(1,4)--cycle;
\draw [fill=pink!40](0,1) -- (1,4)--(1,5)--cycle;
\draw [fill=pink!40](0,1) -- (1,5)--(1,6)--cycle;
\draw [fill=pink!40](0,1) -- (1,6)--(1,7)--cycle;
\draw [fill=pink!40](0,1) -- (1,7)--(1,8)--cycle;
\draw [fill=pink!40](0,1) -- (1,8)--(1,9)--cycle;
\draw [fill=pink!40](0,1) -- (1,9)--(1,10)--cycle;
   \draw [fill=pink!40](0,1) -- (1,10)--(1,11)--cycle;

\draw [fill=pink!40](1,1) -- (1,2)--(2,3)--cycle;
\draw [fill=pink!40](1,1) -- (2,3)--(3,4)--cycle;
\draw [fill=pink!40](1,1) -- (3,4)--(4,5)--cycle;
\draw [fill=pink!40](1,1) -- (4,5)--(5,6)--cycle;
\draw [fill=pink!40](1,1) -- (5,6)--(6,7)--cycle;
\draw [fill=pink!40](1,1) -- (6,7)--(7,8)--cycle;
\draw [fill=pink!40](1,1) -- (7,8)--(8,9)--cycle;
\draw [fill=pink!40](1,1) -- (8,9)--(9,10)--cycle;
   \draw [fill=pink!40](1,1) -- (9,10)--(10,11)--cycle;

\draw [fill=pink!40](1,2) -- (2,3)--(3,5)--cycle;
\draw [fill=pink!40](1,2) -- (3,5)--(4,7)--cycle;
\draw [fill=pink!40](1,2) -- (4,7)--(5,9)--cycle;
   \draw [fill=pink!40](1,2)--(5,9)--(6,11)--cycle;

\draw [fill=yellow](1,2) -- (1,3)--(2,5)--cycle;
\draw [fill=pink!40](1,2) -- (2,5)--(3,7)--cycle;
\draw [fill=pink!40](1,2) -- (3,7)--(4,9)--cycle;
   \draw [fill=pink!40](1,2)--(4,9)--(5,11)--cycle;

\draw [fill=pink!40](2,3) -- (3,5)--(5,8)--cycle;
   \draw [fill=pink!40](2,3)--(5,8)--(7,11)--cycle;

\draw [fill=pink!40](2,3) -- (3,4)--(5,7)--cycle;
\draw [fill=pink!40](2,3) -- (5,7)--(7,10)--cycle;

\draw [fill=pink!40](3,4) -- (4,5)--(7,9)--cycle;

\draw [fill=yellow](1,3) -- (2,5)--(3,8)--cycle;
    \draw [fill=yellow](1,3)--(3,8)--(4,11)--cycle;
\draw [fill=pink!40](1,3) -- (1,4)--(2,7)--cycle;
\draw [fill=pink!40](1,3) -- (2,7)--(3,10)--cycle;

\draw [fill=pink!40](1,4) -- (1,5)--(2,9)--cycle;

\draw [dashed, gray] (0,0) grid (11,11);

\node [above] at (4,11) {$11/4$}; 

\node [below] at (1,0) {$e_{1}$}; 
\node [left] at (0,1) {$e_{2}$}; 
     \node [below] at (4,0) {$4 e_{1}$}; 
     \node [left] at (0, 11) {$11 e_{2}$}; 

 \draw [color=orange!80, line width=1.5pt](0,0) -- (4,11);
  \node[draw,circle,inner sep=1.5pt,fill=black, color=orange] at (4, 11) {};
   \node[draw,circle,inner sep=1.5pt,fill=black, color=black] at (0, 0) {};
   
   \node[draw,circle,inner sep=1.5pt,fill=orange, color=blue] at (1,0) {};
   \node[draw,circle,inner sep=1.5pt,fill=orange, color=red] at (0,1) {};
   \node[draw,circle,inner sep=1.5pt,fill=orange, color=violet!50] at (1,1) {};
   
   \end{scope}


  \begin{scope}[shift={(14,0)}]
  
     \draw [fill=yellow](0,0) -- (6,0)--(3,6)--cycle;
  
      \draw [thick, color=black](0,0) -- (6, 0);
     \draw [thick, color=black](0,0) -- (5.5, 1);
     \draw [thick, color=black](0,0) -- (5,2);
     \draw [thick, color=black](1.5, 3) -- (5,2);
     \draw [thick, color=black](1.5,  3) -- (4.33, 3.33);
     \draw [thick, color=black](1.5, 3) -- (3.66,  4.66);
  
     \draw [thick, color=black](0,0) -- (3,6);
     \draw [thick, color=black](6,0) -- (3,6);
     
     \node [below] at (6,0) {$e_{1}$}; 
     \node [below] at (0,0) {$e_{2}$}; 
     
     \node[draw,circle,inner sep=1.5pt,fill=orange, color=blue] at (6,0) {};
     \node[draw,circle,inner sep=1.5pt,fill=orange, color=red] at (0,0) {};
     
      \node[draw,circle,inner sep=1.5pt,fill=violet!50, color=violet!50] at (5.5,1) {};
     
     \node[draw,circle,inner sep=1.5pt,fill=orange, color=orange] at (3,6) {};

   \end{scope}

   \end{tikzpicture}
\end{center}
  \caption{On the left the lotus $\Lambda(11/4)$ is depicted as a sublotus of the universal lotus 
     and the segment $[0 ,  4 e_1 + 11 e_2]$ giving birth to it. On the right it is depicted as an 
     abstract finite simplicial complex with its base vertices labeled $e_1$ and $e_2$. 
     The color of each vertex of the petal $\delta(e_1, e_2)$ is the color in 
     Figure \ref{fig:embrescusp} of the irreducible component of the total transform of $Z(xy)$ 
     represented by it (see Remark \ref{rem:spacetime}).}  
  \label{fig:Lotus11/4} 
    \end{figure}
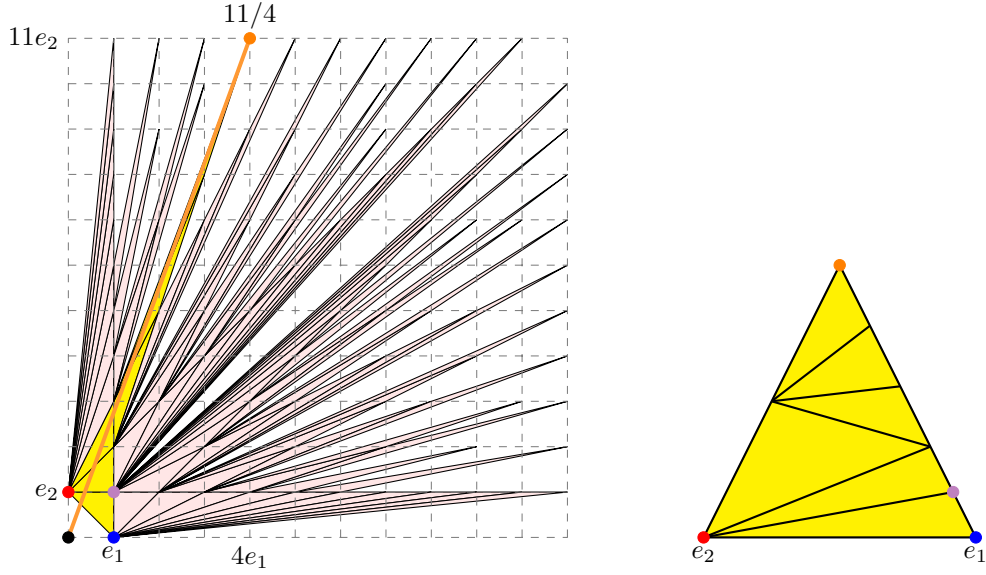

 \begin{figure}[ht!]
     \begin{center}
\begin{tikzpicture}[scale=0.5]

 \begin{scope}[shift={(0,0)}]
 
      \draw [fill=yellow](0,0) -- (6,0)--(3,6)--cycle;
     
     \node[draw,circle,inner sep=1.5pt,fill=orange, color=orange] at (3,6) {};    
     \node [below] at (6,0) {$e_{1}$}; 
     \node [below] at (0,0) {$e_{2}$}; 
      
     \node[draw,circle,inner sep=1.5pt,fill=orange, color=blue] at (6,0) {};
     \node[draw,circle,inner sep=1.5pt,fill=orange, color=red] at (0,0) {};
     
      \draw [->, very thick, color=magenta](6 ,5)   --  (7, 5);  

 \end{scope}

  
   \begin{scope}[shift={(8,0)}]
 
      \draw [fill=yellow](0,0) -- (6,0)--(3,6)--cycle;
  
      \draw [thick, color=black](0,0) -- (6, 0);
     \draw [very thick, color=black](0,0) -- (5,2);
     \draw [very thick, color=black](1.5, 3) -- (5,2);  
     \draw [thick, color=black](0,0) -- (3,6);
     \draw [thick, color=black](6,0) -- (3,6);
     
     \node[draw,circle,inner sep=1.5pt,fill=orange, color=orange] at (3,6) {};
     
     \node [below] at (6,0) {$e_{1}$}; 
     \node [below] at (0,0) {$e_{2}$}; 
     
      \node[draw,circle,inner sep=1.5pt,fill=orange, color=blue] at (6,0) {};
     \node[draw,circle,inner sep=1.5pt,fill=orange, color=red] at (0,0) {};
     
      \draw [->, very thick, color=magenta](6 ,5)   --  (7, 5);  

 \end{scope}


  \begin{scope}[shift={(16,0)}]
  
     \draw [fill=yellow](0,0) -- (6,0)--(3,6)--cycle;
  
      \draw [thick, color=black](0,0) -- (6, 0);
     \draw [thick, color=black](0,0) -- (5.5, 1);
     \draw [very thick, color=black](0,0) -- (5,2);
     \draw [very thick, color=black](1.5, 3) -- (5,2);
     \draw [thick, color=black](1.5,  3) -- (4.33, 3.33);
     \draw [thick, color=black](1.5, 3) -- (3.66,  4.66);
  
     \draw [thick, color=black](0,0) -- (3,6);
     \draw [thick, color=black](6,0) -- (3,6);
     
     \node[draw,circle,inner sep=1.5pt, fill=orange, color=orange] at (3,6) {};
     
     \node [below] at (6,0) {$e_{1}$}; 
     \node [below] at (0,0) {$e_{2}$}; 
     
      \node[draw,circle,inner sep=1.5pt,fill=orange, color=blue] at (6,0) {};
     \node[draw,circle,inner sep=1.5pt,fill=orange, color=red] at (0,0) {};
      
   \end{scope}

   \end{tikzpicture}
\end{center}
  \caption{The construction of the lotus $\Lambda(11/4)$ given the sequence $\mathbf{(2, 1, 3)}$ 
     of quotients of the Euclidean algorithm applied to the pair $(4, 11)$. The composing 
     segments of the zigzag line are thicker than the other segments.}  
  \label{fig:Lotus11/4again} 
    \end{figure}
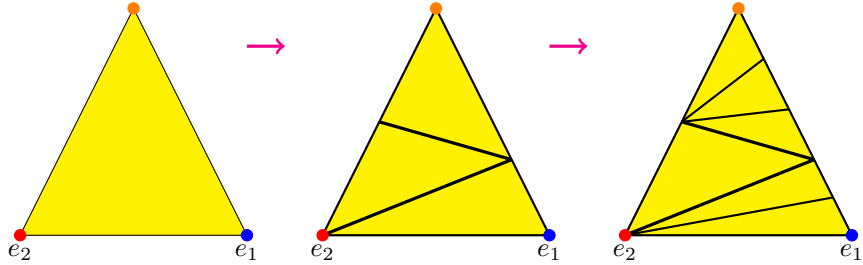    
    
    How to construct the lotus $\Lambda(b/a)$  as a simplicial complex, given the pair $(a,b)$? This is explained in  \cite[Definition 1.5.18, Proposition 1.5.20]{GGP 20}. In Figure \ref{fig:Lotus11/4again} 
    we illustrate this process of construction for the pair $(4,11)$: 
   \begin{itemize} 
       \item[\bf 1)] Apply first the Euclidean algorithm to the pair $(4,11)$, 
           getting the system of equalities:
             \[  \left\{  \begin{array}{cl}
                     11 & = \mathbf{2} \cdot 4 + 3,  \\
                     4  & = \mathbf{1} \cdot 3 + 1, \\
                     3 & = \mathbf{3} \cdot 1.
                 \end{array}   \right.
              \]
       \item[\bf 2)]  Subdivide then into {\bf three} smaller triangles 
            (as many as equalities in the system above) a triangle with two vertices labeled 
             $e_1$ and $e_2$ by a {\bf zigzag line} which starts from the vertex $e_2$. 
       \item[\bf 3)]  Finally, subdivide each of the three triangles in as many smaller triangles 
           as indicated by the corresponding quotient in the Euclidean algorithm. 
           That is, here, into $\mathbf{2, 1, 3}$ smaller triangles. 
    \end{itemize}

     \begin{remark}   \label{rem:initzz}
         The zigzag line appearing in this process is intimately related to the zigzag diagram I used 
          in \cite[Section 5.2]{PP 07} in order to explain geometrically the relation between 
          usual continued fractions and Hirzebruch-Jung ones. Namely, the zigzag line 
          appearing in \cite[Figure 9]{PP 07} is the zigzag line of the lotus obtained 
          by looking of the petals contained in the triangle $[V_0 V_1' A_+]$ of that figure, 
          relative to the basis $(V_0, V_1')$ of the ambient lattice. 
       \end{remark}

Let us extract from formulae  (\ref{eq:seqbu})   the sequence of changes of variables corresponding to the process of embedded resolution of the binomial curve $C_{4, 11}$:

 \begin{equation}   \label{eq:seqbubis}
     \left\{ \begin{array}{l}      
                        x = x_1 \\
                        y = x_1  y_1
                  \end{array} \right.    \rightsquigarrow 
     \left\{ \begin{array}{l}      
                        x_1 = x_2 \\
                        y_1 = x_2  y_2
                  \end{array} \right.      \rightsquigarrow 
     \left\{ \begin{array}{l}      
                        x_2 = x_3 y_3 \\
                        y_2 = y_3
                  \end{array} \right.      \rightsquigarrow    
     \left\{ \begin{array}{l}      
                        x_3 = x_4 \\
                        y_3 = x_4  y_4
                  \end{array} \right.    \rightsquigarrow 
     \left\{ \begin{array}{l}      
                        x_4 = x_5 \\
                        y_4 = x_5  y_5
                  \end{array} \right.        
     \end{equation}
     
     \vspace{2mm}

 \begin{figure}[ht!]
     \begin{center}
\begin{tikzpicture}[scale=0.6]
  
     \draw [fill=yellow](0,0) -- (6,0)--(3,6)--cycle;
  
      \draw [thick, color=black](0,0) -- (6, 0);
     \draw [thick, color=black](0,0) -- (5.5, 1);
     \draw [very thick, color=black](0,0) -- (5,2);
     \draw [very thick, color=black](1.5, 3) -- (5,2);
     \draw [thick, color=black](1.5,  3) -- (4.33, 3.33);
     \draw [thick, color=black](1.5, 3) -- (3.66,  4.66);
  
     \draw [thick, color=black](0,0) -- (3,6);
     \draw [thick, color=black](6,0) -- (3,6);
     
     \node[draw,circle,inner sep=1.5pt,fill=orange, color=orange] at (3,6) {};
     
     \node [below] at (6,0) {$e_{1}$}; 
     \node [below] at (0,0) {$e_{2}$}; 
     
     \node [below] at (4.2, 0.6) {$x$}; 
      \node [below] at (4,1.5) {$x$}; 
     \node [below] at (2,2) {$y$}; 
     \node [below] at (4,3) {$x$}; 
     \node [below] at (3,4) {$x$}; 
      
      \node[draw,circle,inner sep=1.5pt,fill=orange, color=blue] at (6,0) {};
     \node[draw,circle,inner sep=1.5pt,fill=orange, color=red] at (0,0) {};

   \end{tikzpicture}
\end{center}
  \caption{How to read the type of blowup performed at each step, in the sense of Definition \ref{def:affblowup}: excepted for the last petal, write inside each petal either $x$ (if the next petal is attached on its left side) or $y$ (if the next petal is attached on its right side).}  
  \label{fig:Lotus11/4seqxy} 
    \end{figure}
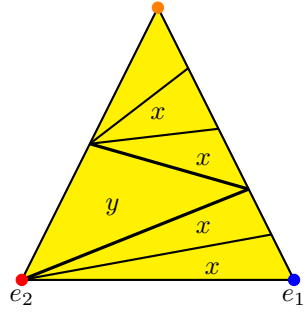

    The previous equalities may be interpreted either as describing successive 
   birational morphisms between affine surfaces, or as true equalities between elements 
   of the fraction fields of those surfaces, identified through the first interpretation. 
   In order to get the minimal embedded resolution of $C_{4, 11}$, one needs 
   to blow up once more, namely, the origin of the plane of coordinates $(x_5, y_5)$.

    We see that the process (\ref{eq:seqbubis}) of blowups consists successively of:
        \begin{itemize}
           \item $2$ $x$-blowups; 
           \item $1$ $y$-blowup;
           \item  $2$ $x$-blowups. 
        \end{itemize}
     This sequence of symbols ``$x$'' and ``$y$'' may be read from the lotus $\Lambda(a/b)$ by starting from the first petal and looking if each new petal is drawn on the left or on the right of the previous one (see Figure \ref{fig:Lotus11/4seqxy}): 
      \begin{itemize}
           \item a left attachment corresponds to an $x$-blowup; 
           \item a right attachment corresponds to an $y$-blowup. 
        \end{itemize}
        
   Note that the last petal of $\Lambda(11/4)$ has no decoration, as it corresponds to 
   the last blowup of the minimal embedded resolution, leading to a strict transform of 
   the curve $C_{4, 11}$ which does not pass through the origins of the two new charts. 
   Therefore, there is no reason to emphasize one of them.

\begin{figure}[h!]
    \begin{center}
\begin{tikzpicture}[scale=1.2]

\begin{scope}[shift={(1,0)}]
        \draw [->, color=red, thick] (-1,0)--(1,0);
        \draw [->, color=blue, thick] (0,-1)--(0,1);

\node[draw,circle, inner sep=1.5pt,color=violet!50, fill=violet!50] at (0,0){};
\node [right, color=blue] at (1, 0) {$x$};
\node [below, color=blue] at (0, -1) {$Z(x)$};
\node [above, color=red] at (0, 1) {$y$};
\node [left, color=red] at (-1, 0) {$Z(y)$};
\end{scope}


\begin{scope}[shift={(8,0)}]
     \draw [-, color=violet!50, very thick](-3,0) -- (3,0);
           \draw [->, color=violet!50, very thick](-2,0) -- (-1,0);
           \draw [->, color=violet!50, very thick](2,0) -- (1,0);
     \draw [->, color=red, thick] (2,-1)--(2,1);
          \node [above, color=violet!50] at (-2, 1) {$\tilde{y}$};
          \node [below, color=red] at (2, -1) {$Z(\bar{y})$};
      \draw [->, color=blue, thick] (-2,-1)--(-2,1);
         \node [above, color=violet!50] at (2, 1) {$\bar{x}$};
         \node [below, color=blue] at (-2, -1) {$Z(\tilde{x})$};
         \node [below, color=blue] at (-1, 0) {$\tilde{x}$};
         \node [below, color=red] at (1, 0) {$\bar{y}$};
         \node [above, color=violet!50] at (0,0) {$E_0$};
\end{scope}


      \draw[<-](3,0)--(4,0);
           \node [above, color=black] at (3.5,0) {$\pi$};
           
\begin{scope}[shift={(-0.5,-2.5)}]
     \draw [-, color=violet!50, very thick](0,0) -- (3,0);
     \node[draw,circle, inner sep=1.5pt,color=blue, fill=blue] at (0,0){};
      \node[draw,circle, inner sep=1.5pt,color=red, fill=red] at (3,0){};
      \node [below, color=blue] at (0, 0) {$Z(x)$};
      \node [below, color=red] at (3, 0) {$Z(y)$};
\end{scope}


\begin{scope}[shift={(6.5,-2.5)}]
     \draw [-, color=violet!50, very thick](0,0) -- (3,0);
     \node[draw,circle, inner sep=1.5pt,color=blue, fill=blue] at (0,0){};
      \node[draw,circle, inner sep=1.5pt,color=red, fill=red] at (3,0){};
      \node[draw,circle, inner sep=1.5pt,color=violet!50, fill=violet!50] at (1.5,0){};
       \node [below, color=blue] at (0, 0) {$Z(x)$};
        \node [below, color=violet!50] at (1.5, 0) {$E_0$};
           \node [above, color=violet!50] at (1.5, 0) {$-1$};
      \node [below, color=red] at (3, 0) {$Z(y)$};
\end{scope}   
           
\end{tikzpicture}
\end{center}
 \caption{The dual graphs of the divisor $Z(xy)$ in $\K^2_{x,y}$ and of its total transform 
    in the blowup surface $\B_0 \K^2$ from Figure \ref{fig:embrescusp}. A way to prove 
    that $E_0 \cdot E_0 = -1$ is explained in Remark \ref{rem:selfint}, point (\ref{it:oncebusi}).}
\label{fig:dgfigure}
   \end{figure}
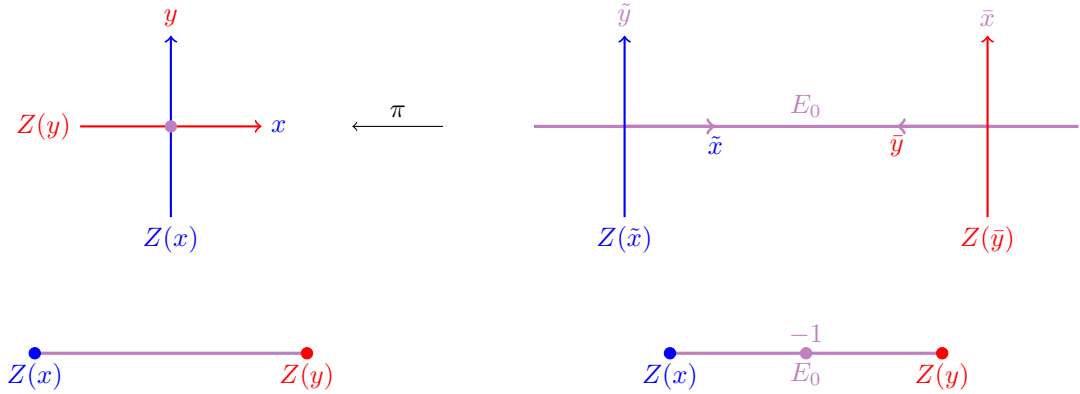

  \medskip
   Let us define now the notion of {\bf dual graph}: 
   
\begin{definition} \label{def:dualgraph}
     Let $D$ be a normal crossings curve contained in a smooth surface $S$, 
     whose irreducible components are smooth. Its {\bf dual graph} is a finite graph 
     constructed as follows (see Figure \ref{fig:dgfigure}):
        \begin{itemize} 
              \item its set of vertices is  in bijection with the set of irreducible components of $D$;  
              \item for every pair of distinct vertices, the set of edges joining them 
                  is in  bijection with the set of their intersection points;
              \item each vertex corresponding to a {\em complete} curve is decorated 
                  by the self-intersection number of that curve in the ambient surface. 
        \end{itemize}
 \end{definition}

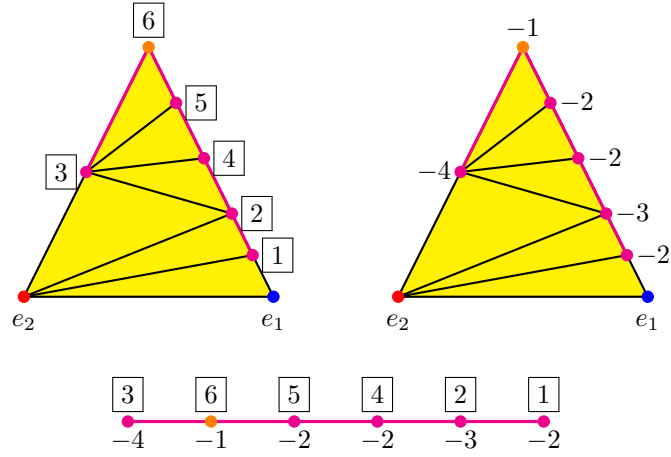
\begin{figure}[ht!]
     \begin{center}
\begin{tikzpicture}[scale=0.55]

   \begin{scope}[shift={(0,0)}]
       \draw [fill=yellow](0,0) -- (6,0)--(3,6)--cycle;
      \draw [thick, color=black](0,0) -- (6, 0);
     \draw [thick, color=black](0,0) -- (5.5, 1);
     \draw [thick, color=black](0,0) -- (5,2);
     \draw [thick, color=black](1.5, 3) -- (5,2);
     \draw [thick, color=black](1.5,  3) -- (4.33, 3.33);
     \draw [thick, color=black](1.5, 3) -- (3.66,  4.66);
  
     \draw [thick, color=black](0,0) -- (3,6);
     \draw [thick, color=black](6,0) -- (3,6);
     
      \draw [very thick, color=magenta](1.5, 3) -- (3,6) -- (5.5, 1);

      \node[draw,circle,inner sep=1.5pt,fill=black, color=blue] at (6, 0) {};
     \node[draw,circle,inner sep=1.5pt,fill=violet, color=red] at (0,0) {};    
     \node[draw,circle,inner sep=1.5pt,fill=black, color=magenta] at (5.5, 1) {};
             \node [right, color=black] at (5.5, 1) {$\boxed{1}$};
     \node[draw,circle,inner sep=1.5pt,fill=red, color=magenta] at (5,2) {};
             \node [right, color=black] at (5, 2) {$\boxed{2}$};
     \node[draw,circle,inner sep=1.5pt,fill=blue, color=magenta] at (1.5, 3) {};
            \node [left, color=black] at (1.5, 3) {$\boxed{3}$};
     \node[draw,circle,inner sep=1.5pt,fill=black, color=magenta] at (4.33, 3.33) {};
              \node [right, color=black] at (4.33, 3.33) {$\boxed{4}$};
     \node[draw,circle,inner sep=1.5pt,fill=black, color=magenta] at (3.66, 4.66) {};
              \node [right, color=black] at (3.66, 4.66) {$\boxed{5}$};
     \node[draw,circle,inner sep=1.5pt,fill=orange, color=orange] at (3,6) {};
              \node [above, color=black] at (3, 6) {$\boxed{6}$};
     
      \node [below, color=black] at (6, -0.2) {$e_1$};
       \node [below, color=black] at (0, -0.2) {$e_2$};

  \end{scope}

  
   \begin{scope}[shift={(9,0)}]
       \draw [fill=yellow](0,0) -- (6,0)--(3,6)--cycle;
      \draw [thick, color=black](0,0) -- (6, 0);
     \draw [thick, color=black](0,0) -- (5.5, 1);
     \draw [thick, color=black](0,0) -- (5,2);
     \draw [thick, color=black](1.5, 3) -- (5,2);
     \draw [thick, color=black](1.5,  3) -- (4.33, 3.33);
     \draw [thick, color=black](1.5, 3) -- (3.66,  4.66);

     \draw [thick, color=black](0,0) -- (3,6);
     \draw [thick, color=black](6,0) -- (3,6);
     
      \draw [very thick, color=magenta](1.5, 3) -- (3,6) -- (5.5, 1);

      \node[draw,circle,inner sep=1.5pt,fill=black, color=blue] at (6, 0) {};
     \node[draw,circle,inner sep=1.5pt,fill=violet, color=red] at (0,0) {};    
     \node[draw,circle,inner sep=1.5pt,fill=black, color=magenta] at (5.5, 1) {};
           \node [right, color=black] at (5.5, 1) {$- 2$};
     \node[draw,circle,inner sep=1.5pt,fill=red, color=magenta] at (5,2) {};
            \node [right, color=black] at (5, 2) {$- 3$};
     \node[draw,circle,inner sep=1.5pt,fill=blue, color=magenta] at (1.5, 3) {};
            \node [left, color=black] at (1.5, 3) {$- 4$};
     \node[draw,circle,inner sep=1.5pt,fill=black, color=magenta] at (4.33, 3.33) {};
             \node [right, color=black] at (4.33, 3.33) {$- 2$};
     \node[draw,circle,inner sep=1.5pt,fill=black, color=magenta] at (3.66, 4.66) {};
             \node [right, color=black] at (3.66, 4.66) {$- 2$};
     \node[draw,circle,inner sep=1.5pt,fill=orange, color=orange] at (3,6) {};
              \node [above, color=black] at (3, 6) {$- 1$};
     
      \node [below, color=black] at (6, -0.2) {$e_1$};
       \node [below, color=black] at (0, -0.2) {$e_2$};

  \end{scope}

   \begin{scope}[shift={(0.5,-3)}]

      \draw [very thick, color=magenta](2, 0) -- (12,0) ;

       \node[draw,circle,inner sep=1.5pt,fill=black, color=magenta] at (2, 0) {};
                 \node [above, color=black] at (2, 0) {$\boxed{3}$};
                 \node [below, color=black] at (2, 0) {$- 4$};
        \node[draw,circle,inner sep=1.5pt,fill=black, color=orange] at (4, 0) {};
               \node [above, color=black] at (4, 0) {$\boxed{6}$};
               \node [below, color=black] at (4, 0) {$- 1$};
         \node[draw,circle,inner sep=1.5pt,fill=black, color=magenta] at (6, 0) {};
               \node [above, color=black] at (6, 0) {$\boxed{5}$};
               \node [below, color=black] at (6, 0) {$- 2$};
          \node[draw,circle,inner sep=1.5pt,fill=black, color=magenta] at (8, 0) {};
                  \node [above, color=black] at (8, 0) {$\boxed{4}$};
                  \node [below, color=black] at (8, 0) {$- 2$};
           \node[draw,circle,inner sep=1.5pt,fill=black, color=magenta] at (10, 0) {};
                   \node [above, color=black] at (10, 0) {$\boxed{2}$};
                   \node [below, color=black] at (10, 0) {$- 3$};
            \node[draw,circle,inner sep=1.5pt,fill=black, color=magenta] at (12, 0) {};
                   \node [above, color=black] at (12, 0) {$\boxed{1}$};
                   \node [below, color=black] at (12, 0) {$- 2$};

  \end{scope}

   \end{tikzpicture}
\end{center}
  \caption{On the upper left is illustrated point {\bf A1)} of Proposition \ref{prop:twoanswers}: 
      the way to see the orders of appearance of the irreducible components of the exceptional 
      divisor $E_{4, 11}$ 
      of the minimal embedded resolution of $C_{4, 11}$.  On the upper right is illustrated point 
      {\bf A2)}: the way to see their self-intersection numbers. On the bottom is drawn the {\em dual 
      graph} of  $E_{4, 11}$ with both types of decorations of its vertices.}  
  \label{fig:Lotus11/4boundary} 
    \end{figure}

       Consider the {\bf dual graph} of the exceptional divisor $\boxed{E_{a,b}}$ 
        of the minimal embedded resolution of a binomial curve $C_{a,b}$. 
       We may answer as follows  the two questions formulated at the end 
       of Section \ref{sect:embres} (see Figure \ref{fig:Lotus11/4boundary}):

\begin{proposition}   \label{prop:twoanswers}  
   Consider the minimal embedded resolution of the binomial curve $C_{a,b}$ and 
   its exceptional divisor $E_{a,b}$.
   \begin{enumerate}
           \item[{\bf A1)}] The dual graph of $E_{a,b}$ 
              is isomorphic to 
              the union of the edges of the boundary of the lotus $\Lambda(a/b)$ which join vertices 
              distinct from $e_1$ and $e_2$. The order of appearance of each 
              irreducible component may be read 
              by looking at the order of attachment of petals, to be marked as a decoration of the 
              corresponding summits. 
          \item[{\bf A2)}] The self-intersection number of the irreducible component of the 
              exceptional divisor $E_{a,b}$ which corresponds to a vertex of the dual graph is equal 
              to the opposite of the number of petals of $\Lambda(a/b)$ 
              which are adjacent to that vertex. 
      \end{enumerate}
 \end{proposition}

\begin{remark}   \label{rem:spacetime}
   Proposition \ref{prop:twoanswers} may be easily proved by induction on $\max(a,b)$, 
   using points (\ref{it:oncebusi}) and (\ref{it:dropsi})  of Remark \ref{rem:selfint}. 
  It holds in fact more generally, for the lotuses associated to any plane curve 
singularity (see \cite[Theorem 1.5.29]{GGP 20} and \cite[Proposition 6.2]{GGP 26}). The vertices 
of such a lotus represent the exceptional divisors created by the successive blowups, as well as 
the coordinate axes of all the charts whose origins were blown up. Two such vertices are 
joined by an edge of the lotus when the corresponding curves intersect in one of the charts. 
Therefore, the lotus may be seen as a {\em space-time representation} of the evolution 
of dual graphs (see \cite[Section 1.7.1, point 10]{GGP 20}). The simplest case of this 
representation may be seen by contemplating the first petal $\delta(e_1, e_2)$ of the lotus 
(see Figure \ref{fig:Lotus11/4}). Its base $[e_1, e_2]$ is the dual graph of the union 
$Z(xy)$ of coordinate axes of the initial affine plane $\K[x,y]$ and the union of its two 
other sides is the dual graph of the total transform of $Z(xy)$ by the blowup morphism $\pi$, 
illustrated on the right side of Figure \ref{fig:dgfigure}. In particular, its summit represents 
the exceptional divisor of the blowup. 
\end{remark}

\begin{remark}  \label{rem:prevdescr}
    This is a sequel of Remark \ref{rem:initzz}.  Before introducing lotuses, in 
    \cite[Section 6.3]{PP 07} I had used the zigzag diagrams of \cite[Section 5.2]{PP 07} 
    in order to describe the dual graphs of the exceptional divisors of the minimal 
    embedded resolutions of the binomial curves $C_{a,b}$. 
\end{remark}

 \medskip
\section{Contracting an anthyphairetic rectangle into a lotus}  \label{sect:contract}
\medskip

In Section \ref{sect:lotposrat} we saw that the lotus $\Lambda(b/a)$ of a positive rational number 
may be constructed from the Euclidean algorithm applied to the pair $(a,b)$. In this section 
we will see an alternative  construction of this lotus, starting from the anthyphairetic process 
applied to the pair $(a,b)$. It passes through an intermediate geometric object, a decorated 
stack of rectangles which we call an {\em anthyphairetic rectangle} $R(a,b)$ 
(see Definition \ref{def:antrect}). 
The lotus $\Lambda(b/a)$ may be obtained by contracting certain  edges of $R(a,b)$ 
(see Proposition \ref{prop:alternconstrlotus}). 
\medskip

  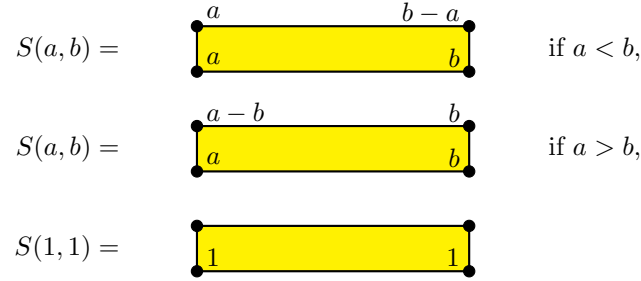
\begin{figure}[h!] 
\begin{center}
\begin{tikzpicture}[x=0.6cm,y=0.6cm]

 \begin{scope}[shift={(0,0)}]
         \node [left, color=black] at (-1.5, 0.5) {$S(a,b) = $};
      \draw [fill=yellow](0,0) -- (6,0)--(6,1)--(0,1)--cycle;
     \draw [thick, color=black](0,0) -- (6,0) -- (6, 1) -- (0,1) -- cycle;
    
     \node[draw,circle,inner sep=1.5pt,fill=black, color=black] at (0,0) {};    
             \node [right, color=black] at (0,0.3) {$a$};
     \node[draw,circle,inner sep=1.5pt,fill=black, color=black] at (0,1) {};
              \node [right, color=black] at (0,1.3) {$a$};
       \node [left, color=black] at (10 ,0.5) {if $a < b$,};

     \node[draw,circle,inner sep=1.5pt,fill=black, color=black] at (6,0) {};
              \node [left, color=black] at (6,0.3) {$b$};
     \node[draw,circle,inner sep=1.5pt,fill=black, color=black] at (6,1) {};
               \node [left, color=black] at (6,1.3) {$b-a$};
     
\end{scope}


  \begin{scope}[shift={(0,-2.2)}]
           \node [left, color=black] at (-1.5, 0.5) {$S(a,b) = $};
       \draw [fill=yellow](0,0) -- (6,0)--(6,1)--(0,1)--cycle;
      \draw [thick, color=black](0,0) -- (6,0) -- (6, 1) -- (0,1) -- cycle;
    
     \node[draw,circle,inner sep=1.5pt,fill=black, color=black] at (0,0) {};    
             \node [right, color=black] at (0,0.3) {$a$};
     \node[draw,circle,inner sep=1.5pt,fill=black, color=black] at (0,1) {};
              \node [right, color=black] at (0,1.3) {$a -b$};
                 
      \node [left, color=black] at (10 ,0.5) {if $a > b$,};

     \node[draw,circle,inner sep=1.5pt,fill=black, color=black] at (6,0) {};
              \node [left, color=black] at (6,0.3) {$b$};
     \node[draw,circle,inner sep=1.5pt,fill=black, color=black] at (6,1) {};
               \node [left, color=black] at (6,1.3) {$b$};
         
   \end{scope}
  
  
   \begin{scope}[shift={(0,-4.4)}]
              \node [left, color=black] at (-1.5, 0.5) {$S(1,1) = $};
       \draw [fill=yellow](0,0) -- (6,0)--(6,1)--(0,1)--cycle;
      \draw [thick, color=black](0,0) -- (6,0) -- (6, 1) -- (0,1) -- cycle;
    
     \node[draw,circle,inner sep=1.5pt,fill=black, color=black] at (0,0) {};    
             \node [right, color=black] at (0,0.3) {$1$};
     \node[draw,circle,inner sep=1.5pt,fill=black, color=black] at (0,1) {};

     \node[draw,circle,inner sep=1.5pt,fill=black, color=black] at (6,0) {};
              \node [left, color=black] at (6,0.3) {$1$};
     \node[draw,circle,inner sep=1.5pt,fill=black, color=black] at (6,1) {};
         
   \end{scope}
  
  
   \end{tikzpicture}
\end{center}
\caption{The three types of subtraction rectangles.}  
   \label{fig:subtrect}
    \end{figure}

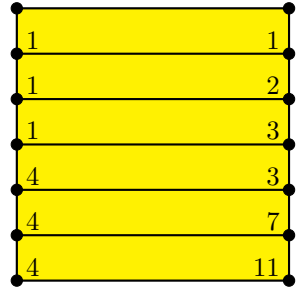
\begin{figure}[h!] 
\begin{center}
\begin{tikzpicture}[x=0.6cm,y=0.6cm]

 \begin{scope}[shift={(0,0)}]
     \draw [fill=yellow](0,0) -- (6,0)--(6,6)--(0,6)--cycle;
     \draw [thick, color=black](0,0) -- (6,0);
     \draw [thick, color=black](0,1) -- (6,1);
     \draw [thick, color=black](0,2) -- (6,2);
     \draw [thick, color=black](0,3) -- (6,3);
     \draw [thick, color=black](0,4) -- (6,4);
     \draw [thick, color=black](0,5) -- (6,5);
     \draw [thick, color=black](0,6) -- (6,6);

     \draw [thick, color=black](0,0) -- (0,6);
     \draw [thick, color=black](6,0) -- (6,6);

     \node[draw,circle,inner sep=1.5pt,fill=black, color=black] at (0,0) {};    
             \node [right, color=black] at (0,0.3) {$4$};
     \node[draw,circle,inner sep=1.5pt,fill=black, color=black] at (0,1) {};
              \node [right, color=black] at (0,1.3) {$4$};
     \node[draw,circle,inner sep=1.5pt,fill=black, color=black] at (0,2) {};
               \node [right, color=black] at (0,2.3) {$4$};
     \node[draw,circle,inner sep=1.5pt,fill=black, color=black] at (0,3) {};
                \node [right, color=black] at (0,3.3) {$1$};
     \node[draw,circle,inner sep=1.5pt,fill=black, color=black] at (0,4) {};
               \node [right, color=black] at (0,4.3) {$1$};
     \node[draw,circle,inner sep=1.5pt,fill=black, color=black] at (0,5) {};
               \node [right, color=black] at (0,5.3) {$1$};
     \node[draw,circle,inner sep=1.5pt,fill=black, color=black] at (0,6) {};

     \node[draw,circle,inner sep=1.5pt,fill=black, color=black] at (6,0) {};
              \node [left, color=black] at (6,0.3) {$11$};
     \node[draw,circle,inner sep=1.5pt,fill=black, color=black] at (6,1) {};
               \node [left, color=black] at (6,1.3) {$7$};
     \node[draw,circle,inner sep=1.5pt,fill=black, color=black] at (6,2) {};
                \node [left, color=black] at (6,2.3) {$3$};
     \node[draw,circle,inner sep=1.5pt,fill=black, color=black] at (6,3) {};
                \node [left, color=black] at (6,3.3) {$3$};
     \node[draw,circle,inner sep=1.5pt,fill=black, color=black] at (6,4) {};
                \node [left, color=black] at (6,4.3) {$2$};
     \node[draw,circle,inner sep=1.5pt,fill=black, color=black] at (6,5) {};
                 \node [left, color=black] at (6,5.3) {$1$};
     \node[draw,circle,inner sep=1.5pt,fill=black, color=black] at (6,6) {};
\end{scope}

     
       \end{tikzpicture}
\end{center}
\caption{The anthyphairetic rectangle $R(4, 11)$ of the pair $(4,11)$ 
    (see Definition \ref{def:antrect}).}  
   \label{fig:anthyrect}
    \end{figure}

       \begin{definition}   \label{def:antrect}
        Let $(a,b)$ be a pair of coprime positive integers.  The {\bf subtraction rectangle} 
        $\boxed{S(a,b)}$ of $(a,b)$  is defined as shown in Figure \ref{fig:subtrect}.               
        The {\bf anthyphairetic rectangle $\boxed{R(a,b)}$ of $(a,b)$} is obtained 
        as shown in Figure \ref{fig:anthyrect}, by stacking the subtraction rectangles 
        of all the pairs of coprime positive integers obtained by applying the anthyphairetic process 
        to $(a,b)$. 
   \end{definition}

\begin{figure}[ht!] 
\begin{center}
\begin{tikzpicture}[x=0.6cm,y=0.6cm]

 \begin{scope}[shift={(0,0)}]
     \draw [fill=yellow](0,0) -- (6,0)--(6,6)--(0,6)--cycle;
     \draw [thick, color=black](0,0) -- (6,0);
     \draw [thick, color=black](0,1) -- (6,1);
     \draw [thick, color=black](0,2) -- (6,2);
     \draw [thick, color=black](0,3) -- (6,3);
     \draw [thick, color=black](0,4) -- (6,4);
     \draw [thick, color=black](0,5) -- (6,5);
     \draw [thick, color=black](0,6) -- (6,6);

     \draw [thick, color=black](0,0) -- (0,6);
     \draw [thick, color=black](6,0) -- (6,6);

     \node[draw,circle,inner sep=1.5pt,fill=black, color=black] at (0,0) {};    
             \node [right, color=black] at (0,0.3) {$4$};
     \node[draw,circle,inner sep=1.5pt,fill=black, color=black] at (0,1) {};
              \node [right, color=black] at (0,1.3) {$4$};
     \node[draw,circle,inner sep=1.5pt,fill=black, color=black] at (0,2) {};
               \node [right, color=black] at (0,2.3) {$4$};
     \node[draw,circle,inner sep=1.5pt,fill=black, color=black] at (0,3) {};
                \node [right, color=black] at (0,3.3) {$1$};
     \node[draw,circle,inner sep=1.5pt,fill=black, color=black] at (0,4) {};
               \node [right, color=black] at (0,4.3) {$1$};
     \node[draw,circle,inner sep=1.5pt,fill=black, color=black] at (0,5) {};
               \node [right, color=black] at (0,5.3) {$1$};
     \node[draw,circle,inner sep=1.5pt,fill=black, color=black] at (0,6) {};

     \node[draw,circle,inner sep=1.5pt,fill=black, color=black] at (6,0) {};
              \node [left, color=black] at (6,0.3) {$11$};
     \node[draw,circle,inner sep=1.5pt,fill=black, color=black] at (6,1) {};
               \node [left, color=black] at (6,1.3) {$7$};
     \node[draw,circle,inner sep=1.5pt,fill=black, color=black] at (6,2) {};
                \node [left, color=black] at (6,2.3) {$3$};
     \node[draw,circle,inner sep=1.5pt,fill=black, color=black] at (6,3) {};
                \node [left, color=black] at (6,3.3) {$3$};
     \node[draw,circle,inner sep=1.5pt,fill=black, color=black] at (6,4) {};
                \node [left, color=black] at (6,4.3) {$2$};
     \node[draw,circle,inner sep=1.5pt,fill=black, color=black] at (6,5) {};
                 \node [left, color=black] at (6,5.3) {$1$};
     \node[draw,circle,inner sep=1.5pt,fill=black, color=black] at (6,6) {};
     
          \draw [->, very thick, color=magenta](7.2 ,4)   --  (8.2, 4);  
     
\end{scope}


  \begin{scope}[shift={(9.5,0)}]
      \draw [fill=yellow](0,0) -- (6,0)--(6,6)--(0,6)--cycle;
     \draw [thick, color=black](0,0) -- (6,0);
     \draw [thick, color=black](0,1) -- (6,1);
     \draw [thick, color=black](0,2) -- (6,2);
     \draw [thick, color=black](0,3) -- (6,3);
     \draw [thick, color=black](0,4) -- (6,4);
     \draw [thick, color=black](0,5) -- (6,5);
     \draw [thick, color=orange, line width=2.5pt](0,6) -- (6,6);

     \draw [thick, color=black](0,0) -- (0,6);
     \draw [thick, color=black](6,0) -- (6,6);
     
      \draw [color=violet, line width=2.5pt](0,0) -- (0,2);
       \draw [color=cyan, line width=2.5pt](0,3) -- (0,5);
      
        \draw [color=magenta, line width=2.5pt](6,2) -- (6,3);

     \node[draw,circle,inner sep=1.5pt,fill=black, color=black] at (0,0) {};    
             \node [right, color=violet] at (0,0.3) {$4$};
     \node[draw,circle,inner sep=1.5pt,fill=black, color=black] at (0,1) {};
              \node [right, color=violet] at (0,1.3) {$4$};
     \node[draw,circle,inner sep=1.5pt,fill=black, color=black] at (0,2) {};
               \node [right, color=violet] at (0,2.3) {$4$};
     \node[draw,circle,inner sep=1.5pt,fill=black, color=black] at (0,3) {};
                \node [right, color=cyan] at (0,3.3) {$1$};
     \node[draw,circle,inner sep=1.5pt,fill=black, color=black] at (0,4) {};
               \node [right, color=cyan] at (0,4.3) {$1$};
     \node[draw,circle,inner sep=1.5pt,fill=black, color=black] at (0,5) {};
               \node [right, color=cyan] at (0,5.3) {$1$};
     \node[draw,circle,inner sep=1.5pt,fill=black, color=black] at (0,6) {};

     \node[draw,circle,inner sep=1.5pt,fill=black, color=black] at (6,0) {};
              \node [left, color=black] at (6,0.3) {$11$};
     \node[draw,circle,inner sep=1.5pt,fill=black, color=black] at (6,1) {};
               \node [left, color=black] at (6,1.3) {$7$};
     \node[draw,circle,inner sep=1.5pt,fill=black, color=black] at (6,2) {};
                \node [left, color=magenta] at (6,2.3) {$3$};
     \node[draw,circle,inner sep=1.5pt,fill=black, color=black] at (6,3) {};
                \node [left, color=magenta] at (6,3.3) {$3$};
     \node[draw,circle,inner sep=1.5pt,fill=black, color=black] at (6,4) {};
                \node [left, color=black] at (6,4.3) {$2$};
     \node[draw,circle,inner sep=1.5pt,fill=black, color=black] at (6,5) {};
                 \node [left, color=black] at (6,5.3) {$1$};
     \node[draw,circle,inner sep=1.5pt,fill=black, color=black] at (6,6) {};
     
            \draw [->, very thick, color=magenta](7.5 ,4)   --  (8.5, 4); 
  \end{scope}
  
  
   \begin{scope}[shift={(18.5,0)}]
       \draw [fill=yellow](0,0) -- (6,0)--(3,6)--cycle;
      \draw [thick, color=black](0,0) -- (6, 0);
     \draw [thick, color=black](0,0) -- (5.5, 1);
     \draw [thick, color=black](0,0) -- (5,2);
     \draw [thick, color=black](1.5, 3) -- (5,2);
     \draw [thick, color=black](1.5,  3) -- (4.33, 3.33);
     \draw [thick, color=black](1.5, 3) -- (3.66,  4.66);
  
     \draw [thick, color=black](0,0) -- (3,6);
     \draw [thick, color=black](6,0) -- (3,6);
    
      \node[draw,circle,inner sep=1.5pt,fill=black, color=black] at (6, 0) {};
                  \node [left, color=black] at (6, 0.3) {$11$};
     \node[draw,circle,inner sep=1.5pt,fill=violet, color=violet] at (0,0) {};    
               \node [right, color=violet] at (0.2,0.5) {$4$};
     \node[draw,circle,inner sep=1.5pt,fill=black, color=black] at (5.5, 1) {};
                  \node [left, color=black] at (5.4, 1.3) {$7$};
     \node[draw,circle,inner sep=1.5pt,fill=red, color=magenta] at (5,2) {};
                 \node [left, color=magenta] at (4.5,2) {$3$};
     \node[draw,circle,inner sep=1.5pt,fill=blue, color=cyan] at (1.5, 3) {};
                 \node [right, color=cyan] at (1.5, 2.6) {$1$};
     \node[draw,circle,inner sep=1.5pt,fill=black, color=black] at (4.33, 3.33) {};
                  \node [left, color=black] at (4.3, 3.6) {$2$};
     \node[draw,circle,inner sep=1.5pt,fill=black, color=black] at (3.66, 4.66) {};
                  \node [left, color=black] at (3.66, 4.66) {$1$};
     \node[draw,circle,inner sep=1.5pt,fill=orange, color=orange] at (3,6) {};
     
      \node [below, color=black] at (6, -0.2) {$e_1$};
       \node [below, color=black] at (0, -0.2) {$e_2$};

  \end{scope}

 \end{tikzpicture}
\end{center}
\caption{Contracting the anthyphairetic rectangle $R(4,11)$ of Figure \ref{fig:anthyrect} 
     into the anthyphairetic triangle $T(4,11)$, as explained in Definition \ref{def:triangpair}.}  
   \label{fig:contrant}
    \end{figure}
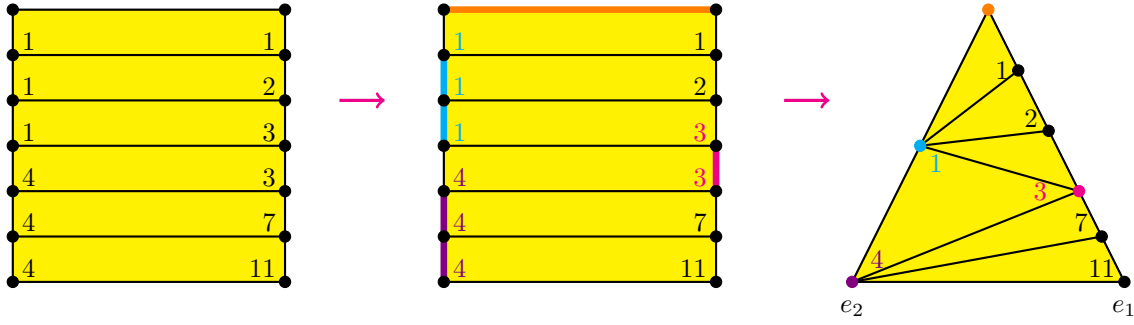

    The positive integers decorating the vertices of the anthyphairetic rectangle allow to exhibit a finite set of  edges which have to be contracted:

 \begin{definition}   \label{def:triangpair}
        Let $(a,b)$ be a pair of coprime positive integers.  The {\bf anthyphairetic triangle} 
        $\boxed{T(a,b)}$ of $(a,b)$ is obtained from the anthyphairetic rectangle 
        $R(a,b)$ of $(a,b)$ by contracting the following edges to points (see Figure \ref{fig:contrant}):
           \begin{itemize}
               \item the vertical edges whose vertices have the same numerical decoration;
               \item the top edge of the subtraction rectangle $S(1,1)$. 
           \end{itemize} 
        The lower right vertex is decorated by $e_1$ and the lower left one by $e_2$. One 
        decorates each vertex of $T(a,b)$ by the positive integer which decorates any 
        of its preimage vertices in $R(a,b)$. 
   \end{definition}

 One may easily prove by induction on $\max(a,b)$ that one gets in this way an alternative 
 construction of the lotus of the positive rational number $b/a$:
 
 \begin{proposition}  \label{prop:alternconstrlotus}
      Let $(a,b)$ be a pair of coprime positive integers.  The anthyphairetic triangle 
        $T(a,b)$ of $(a,b)$ is isomorphic as an undecorated  simplicial complex to the lotus 
        $\Lambda(b/a)$ of the rational number $b/a$. 
 \end{proposition}

  \begin{figure}[ht!] 
\begin{center}
\begin{tikzpicture}[scale=0.5] 

 \begin{scope}[shift={(0,0)}]
      \draw [fill=yellow](0,0) -- (6,0)--(6,6)--(0,6)--cycle;
     \draw [thick, color=black](0,0) -- (6,0);
     \draw [thick, color=black](0,1) -- (6,1);
     \draw [thick, color=black](0,2) -- (6,2);
     \draw [thick, color=black](0,3) -- (6,3);
     \draw [thick, color=black](0,4) -- (6,4);
     \draw [thick, color=black](0,5) -- (6,5);
     \draw [thick, color=black](0,6) -- (6,6);

     \draw [thick, color=black](0,0) -- (0,6);
     \draw [thick, color=black](6,0) -- (6,6);

     \node[draw,circle,inner sep=1.5pt,fill=black, color=black] at (0,0) {};    
             \node [right, color=black] at (0,0.3) {$4$};
               \node [left, color=black] at (0,0) {$x$};
     \node[draw,circle,inner sep=1.5pt,fill=black, color=black] at (0,1) {};
              \node [right, color=black] at (0,1.3) {$4$};
                 \node [left, color=black] at (0,1) {$x_1$};
     \node[draw,circle,inner sep=1.5pt,fill=black, color=black] at (0,2) {};
               \node [right, color=black] at (0,2.3) {$4$};
                  \node [left, color=black] at (0,2) {$x_2$};
     \node[draw,circle,inner sep=1.5pt,fill=black, color=black] at (0,3) {};
                \node [right, color=black] at (0,3.3) {$1$};
                 \node [left, color=black] at (0,3) {$x_3$};
     \node[draw,circle,inner sep=1.5pt,fill=black, color=black] at (0,4) {};
               \node [right, color=black] at (0,4.3) {$1$};
                  \node [left, color=black] at (0,4) {$x_4$};
     \node[draw,circle,inner sep=1.5pt,fill=black, color=black] at (0,5) {};
               \node [right, color=black] at (0,5.3) {$1$};
                  \node [left, color=black] at (0,5) {$x_5$};
     \node[draw,circle,inner sep=1.5pt,fill=black, color=black] at (0,6) {};

     \node[draw,circle,inner sep=1.5pt,fill=black, color=black] at (6,0) {};
              \node [left, color=black] at (6,0.3) {$11$};
                  \node [right, color=black] at (6,0) {$y$};
     \node[draw,circle,inner sep=1.5pt,fill=black, color=black] at (6,1) {};
               \node [left, color=black] at (6,1.3) {$7$};
                 \node [right, color=black] at (6,1) {$y_1$};
     \node[draw,circle,inner sep=1.5pt,fill=black, color=black] at (6,2) {};
                \node [left, color=black] at (6,2.3) {$3$};
                  \node [right, color=black] at (6,2) {$y_2$};
     \node[draw,circle,inner sep=1.5pt,fill=black, color=black] at (6,3) {};
                \node [left, color=black] at (6,3.3) {$3$};
                  \node [right, color=black] at (6,3) {$y_3$};
     \node[draw,circle,inner sep=1.5pt,fill=black, color=black] at (6,4) {};
                \node [left, color=black] at (6,4.3) {$2$};
                  \node [right, color=black] at (6,4) {$y_4$};
     \node[draw,circle,inner sep=1.5pt,fill=black, color=black] at (6,5) {};
                 \node [left, color=black] at (6,5.3) {$1$};
                   \node [right, color=black] at (6,5) {$y_5$};
     \node[draw,circle,inner sep=1.5pt,fill=black, color=black] at (6,6) {};
     
     \draw [->, very thick, color=magenta](7.5 ,4.5)   --  (8.5, 4.5); 
\end{scope}


  \begin{scope}[shift={(10,0)}]
      \draw [fill=yellow](0,0) -- (6,0)--(6,6)--(0,6)--cycle;
     \draw [thick, color=black](0,0) -- (6,0);
     \draw [thick, color=black](0,1) -- (6,1);
     \draw [thick, color=black](0,2) -- (6,2);
     \draw [thick, color=black](0,3) -- (6,3);
     \draw [thick, color=black](0,4) -- (6,4);
     \draw [thick, color=black](0,5) -- (6,5);
     \draw [thick, color=orange, line width=2.5pt](0,6) -- (6,6);

     \draw [thick, color=black](0,0) -- (0,6);
     \draw [thick, color=black](6,0) -- (6,6);
     
      \draw [color=violet, line width=2.5pt](0,0) -- (0,2);
       \draw [color=cyan, line width=2.5pt](0,3) -- (0,5);
      
        \draw [color=magenta, line width=2.5pt](6,2) -- (6,3);

     \node[draw,circle,inner sep=1.5pt,fill=black, color=black] at (0,0) {};    
             \node [right, color=violet] at (0,0.3) {$4$};
               \node [left, color=black] at (0,0) {$x$};
     \node[draw,circle,inner sep=1.5pt,fill=black, color=black] at (0,1) {};
              \node [right, color=violet] at (0,1.3) {$4$};
                 \node [left, color=black] at (0,1) {$x_1$};
     \node[draw,circle,inner sep=1.5pt,fill=black, color=black] at (0,2) {};
               \node [right, color=violet] at (0,2.3) {$4$};
                  \node [left, color=black] at (0,2) {$x_2$};
     \node[draw,circle,inner sep=1.5pt,fill=black, color=black] at (0,3) {};
                \node [right, color=cyan] at (0,3.3) {$1$};
                 \node [left, color=black] at (0,3) {$x_3$};
     \node[draw,circle,inner sep=1.5pt,fill=black, color=black] at (0,4) {};
               \node [right, color=cyan] at (0,4.3) {$1$};
                  \node [left, color=black] at (0,4) {$x_4$};
     \node[draw,circle,inner sep=1.5pt,fill=black, color=black] at (0,5) {};
               \node [right, color=cyan] at (0,5.3) {$1$};
                  \node [left, color=black] at (0,5) {$x_5$};
     \node[draw,circle,inner sep=1.5pt,fill=black, color=black] at (0,6) {};

     \node[draw,circle,inner sep=1.5pt,fill=black, color=black] at (6,0) {};
              \node [left, color=black] at (6,0.3) {$11$};
                  \node [right, color=black] at (6,0) {$y$};
     \node[draw,circle,inner sep=1.5pt,fill=black, color=black] at (6,1) {};
               \node [left, color=black] at (6,1.3) {$7$};
                 \node [right, color=black] at (6,1) {$y_1$};
     \node[draw,circle,inner sep=1.5pt,fill=black, color=black] at (6,2) {};
                \node [left, color=magenta] at (6,2.3) {$3$};
                  \node [right, color=black] at (6,2) {$y_2$};
     \node[draw,circle,inner sep=1.5pt,fill=black, color=black] at (6,3) {};
                \node [left, color=magenta] at (6,3.3) {$3$};
                  \node [right, color=black] at (6,3) {$y_3$};
     \node[draw,circle,inner sep=1.5pt,fill=black, color=black] at (6,4) {};
                \node [left, color=black] at (6,4.3) {$2$};
                  \node [right, color=black] at (6,4) {$y_4$};
     \node[draw,circle,inner sep=1.5pt,fill=black, color=black] at (6,5) {};
                 \node [left, color=black] at (6,5.3) {$1$};
                   \node [right, color=black] at (6,5) {$y_5$};
     \node[draw,circle,inner sep=1.5pt,fill=black, color=black] at (6,6) {};
     
      \draw [->, very thick, color=magenta](8.5 ,4.5)   --  (9.5, 4.5); 
  \end{scope}
  
  
   \begin{scope}[shift={(22,0)}]
        \draw [fill=yellow](0,0) -- (6,0)--(3,6)--cycle;
      \draw [thick, color=black](0,0) -- (6, 0);
     \draw [thick, color=black](0,0) -- (5.5, 1);
     \draw [thick, color=black](0,0) -- (5,2);
     \draw [thick, color=black](1.5, 3) -- (5,2);
     \draw [thick, color=black](1.5,  3) -- (4.33, 3.33);
     \draw [thick, color=black](1.5, 3) -- (3.66,  4.66);
  
     \draw [thick, color=black](0,0) -- (3,6);
     \draw [thick, color=black](6,0) -- (3,6);
    
      \node[draw,circle,inner sep=1.5pt,fill=black, color=black] at (6, 0) {};
                  \node [right, color=black] at (6, 0) {$y$};
     \node[draw,circle,inner sep=1.5pt,fill=violet, color=violet] at (0,0) {};    
               \node [left, color=violet] at (0,0) {$x = x_1 = x_2$};
     \node[draw,circle,inner sep=1.5pt,fill=black, color=black] at (5.5, 1) {};
                  \node [right, color=black] at (5.5, 1) {$y_1$};
     \node[draw,circle,inner sep=1.5pt,fill=red, color=magenta] at (5,2) {};
                 \node [right, color=magenta] at (5,2) {$y_2 = y_3$};
     \node[draw,circle,inner sep=1.5pt,fill=blue, color=cyan] at (1.5, 3) {};
                 \node [left, color=cyan] at (1.5, 3) {$x_3= x_4 = x_5$};
     \node[draw,circle,inner sep=1.5pt,fill=black, color=black] at (4.33, 3.33) {};
                  \node [right, color=black] at (4.33, 3.33) {$y_4$};
     \node[draw,circle,inner sep=1.5pt,fill=black, color=black] at (3.66, 4.66) {};
                  \node [right, color=black] at (3.66, 4.66) {$y_5$};
     \node[draw,circle,inner sep=1.5pt,fill=orange, color=orange] at (3,6) {};

  \end{scope}

 \end{tikzpicture}
\end{center}
\caption{Decorating the anthyphairetic rectangle $R(4,11)$ and its contracted 
        anthyphairetic triangle $T(4,11)$ using the coordinates of the charts appeared during 
        the blowup process of the binomial curve $C_{4, 11}$ (see Remark \ref{rem:decxiyi}).}  
   \label{fig:coordec}
    \end{figure}
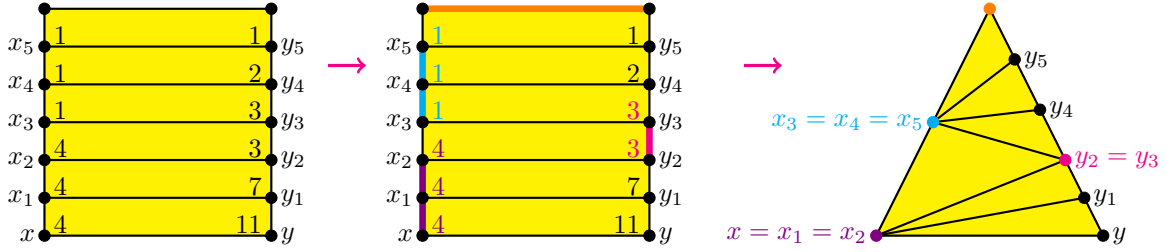

 \begin{remark}  \label{rem:decxiyi}
        Note that if one decorates the left side of the anthyphairetic rectangle $R(4, 11)$ 
        by the $x$-variables 
        $x, x_1, x_2, \dots$ and the right side by the $y$-variables $y, y_1, y_2, \dots$ 
        as illustrated on the left of Figure \ref{fig:coordec}, then one gets as decorations 
        of the vertices which are identified by the contraction of Definition \ref{def:triangpair} 
        exactly the variables which are identified by the changes of variables (\ref{eq:seqbubis}) 
        of the blowup process. This fact is again completely general, and 
        may be also proved by induction on $\max(a,b)$. 
   \end{remark}

\section{Transforming a lotus into an anthyphairetic rectangle}  \label{sec:lotintorect}

In the previous section we have seen how to pass graphically from the anthyphairetic process applied 
 to a coprime pair $(a,b)$ to the corresponding continued fraction, by transforming the 
 anthyphairetic rectangle $R(a,b)$ into the anthyphairetic triangle $T(a,b)$, 
 which is isomorphic to the lotus 
 $\Lambda(b/a)$. In this section we will examine a converse process. 
 \medskip

 \begin{figure}[ht!] 
\begin{center}
\begin{tikzpicture}[scale=0.55]

   \begin{scope}[shift={(0,0)}]
       \draw [fill=yellow](0,0) -- (6,0)--(3,6)--cycle;
      \draw [thick, color=black](0,0) -- (6, 0);
     \draw [thick, color=black](0,0) -- (5.5, 1);
     \draw [thick, color=black](0,0) -- (5,2);
     \draw [thick, color=black](1.5, 3) -- (5,2);
     \draw [thick, color=black](1.5,  3) -- (4.33, 3.33);
     \draw [thick, color=black](1.5, 3) -- (3.66,  4.66);
  
     \draw [thick, color=black](0,0) -- (3,6);
     \draw [thick, color=black](6,0) -- (3,6);

      \node[draw,circle,inner sep=1.5pt,fill=black, color=black] at (6, 0) {};
     \node[draw,circle,inner sep=1.5pt,fill=violet, color=black] at (0,0) {};    
     \node[draw,circle,inner sep=1.5pt,fill=black, color=black] at (5.5, 1) {};
     \node[draw,circle,inner sep=1.5pt,fill=red, color=black] at (5,2) {};
     \node[draw,circle,inner sep=1.5pt,fill=blue, color=black] at (1.5, 3) {};
                 \node [right, color=black] at (1.5, 2.6) {$1$};
     \node[draw,circle,inner sep=1.5pt,fill=black, color=black] at (4.33, 3.33) {};
     \node[draw,circle,inner sep=1.5pt,fill=black, color=black] at (3.66, 4.66) {};
                  \node [left, color=black] at (3.66, 4.66) {$1$};
     \node[draw,circle,inner sep=1.5pt,fill=orange, color=orange] at (3,6) {};
     
      \node [below, color=black] at (6, -0.2) {$e_1$};
       \node [below, color=black] at (0, -0.2) {$e_2$};
       
        \draw [->, very thick, color=magenta](7 ,5)   --  (8, 5);  

  \end{scope}

  
   \begin{scope}[shift={(9,0)}]
       \draw [fill=yellow](0,0) -- (6,0)--(3,6)--cycle;
      \draw [thick, color=black](0,0) -- (6, 0);
     \draw [thick, color=black](0,0) -- (5.5, 1);
     \draw [thick, color=black](0,0) -- (5,2);
     \draw [thick, color=black](1.5, 3) -- (5,2);
     \draw [thick, color=black](1.5,  3) -- (4.33, 3.33);
     \draw [thick, color=black](1.5, 3) -- (3.66,  4.66);
  
     \draw [thick, color=black](0,0) -- (3,6);
     \draw [thick, color=black](6,0) -- (3,6);

      \node[draw,circle,inner sep=1.5pt,fill=black, color=black] at (6, 0) {};
     \node[draw,circle,inner sep=1.5pt,fill=violet, color=black] at (0,0) {};    
     \node[draw,circle,inner sep=1.5pt,fill=black, color=black] at (5.5, 1) {};
     \node[draw,circle,inner sep=1.5pt,fill=red, color=black] at (5,2) {};
     \node[draw,circle,inner sep=1.5pt,fill=blue, color=black] at (1.5, 3) {};
                 \node [right, color=black] at (1.5, 2.6) {$1$};
     \node[draw,circle,inner sep=1.5pt,fill=black, color=cyan] at (4.33, 3.33) {};
                  \node [left, color=cyan] at (4.3, 3.6) {$2$};
     \node[draw,circle,inner sep=1.5pt,fill=black, color=black] at (3.66, 4.66) {};
                  \node [left, color=black] at (3.66, 4.66) {$1$};
     \node[draw,circle,inner sep=1.5pt,fill=orange, color=orange] at (3,6) {};
     
      \node [below, color=black] at (6, -0.2) {$e_1$};
       \node [below, color=black] at (0, -0.2) {$e_2$};
       
        \draw [->, very thick, color=magenta](7 ,5)   --  (8, 5);  

  \end{scope}
  
  
   \begin{scope}[shift={(18.5,0)}]
       \draw [fill=yellow](0,0) -- (6,0)--(3,6)--cycle;
      \draw [thick, color=black](0,0) -- (6, 0);
     \draw [thick, color=black](0,0) -- (5.5, 1);
     \draw [thick, color=black](0,0) -- (5,2);
     \draw [thick, color=black](1.5, 3) -- (5,2);
     \draw [thick, color=black](1.5,  3) -- (4.33, 3.33);
     \draw [thick, color=black](1.5, 3) -- (3.66,  4.66);
  
     \draw [thick, color=black](0,0) -- (3,6);
     \draw [thick, color=black](6,0) -- (3,6);
    
      \node[draw,circle,inner sep=1.5pt,fill=black, color=black] at (6, 0) {};
     \node[draw,circle,inner sep=1.5pt,fill=violet, color=black] at (0,0) {};    
     \node[draw,circle,inner sep=1.5pt,fill=black, color=black] at (5.5, 1) {};
     \node[draw,circle,inner sep=1.5pt,fill=red, color=cyan] at (5,2) {};
                 \node [left, color=cyan] at (4.5,2) {$3$};
     \node[draw,circle,inner sep=1.5pt,fill=blue, color=black] at (1.5, 3) {};
                 \node [right, color=black] at (1.5, 2.6) {$1$};
     \node[draw,circle,inner sep=1.5pt,fill=black, color=black] at (4.33, 3.33) {};
                  \node [left, color=black] at (4.3, 3.6) {$2$};
     \node[draw,circle,inner sep=1.5pt,fill=black, color=black] at (3.66, 4.66) {};
                  \node [left, color=black] at (3.66, 4.66) {$1$};
     \node[draw,circle,inner sep=1.5pt,fill=orange, color=orange] at (3,6) {};
     
      \node [below, color=black] at (6, -0.2) {$e_1$};
       \node [below, color=black] at (0, -0.2) {$e_2$};
       
        \draw [-> , very thick, color=magenta](3 , -1.5)   --  (3, -2.5);  

  \end{scope}
  
  
   \begin{scope}[shift={(18.5, -9.5)}]
       \draw [fill=yellow](0,0) -- (6,0)--(3,6)--cycle;
      \draw [thick, color=black](0,0) -- (6, 0);
     \draw [thick, color=black](0,0) -- (5.5, 1);
     \draw [thick, color=black](0,0) -- (5,2);
     \draw [thick, color=black](1.5, 3) -- (5,2);
     \draw [thick, color=black](1.5,  3) -- (4.33, 3.33);
     \draw [thick, color=black](1.5, 3) -- (3.66,  4.66);
  
     \draw [thick, color=black](0,0) -- (3,6);
     \draw [thick, color=black](6,0) -- (3,6);
    
      \node[draw,circle,inner sep=1.5pt,fill=black, color=black] at (6, 0) {};
     \node[draw,circle,inner sep=1.5pt,fill=violet, color=cyan] at (0,0) {};    
               \node [right, color=cyan] at (0.2,0.5) {$4$};
     \node[draw,circle,inner sep=1.5pt,fill=black, color=black] at (5.5, 1) {};
     \node[draw,circle,inner sep=1.5pt,fill=red, color=black] at (5,2) {};
                 \node [left, color=black] at (4.5,2) {$3$};
     \node[draw,circle,inner sep=1.5pt,fill=blue, color=black] at (1.5, 3) {};
                 \node [right, color=black] at (1.5, 2.6) {$1$};
     \node[draw,circle,inner sep=1.5pt,fill=black, color=black] at (4.33, 3.33) {};
                  \node [left, color=black] at (4.3, 3.6) {$2$};
     \node[draw,circle,inner sep=1.5pt,fill=black, color=black] at (3.66, 4.66) {};
                  \node [left, color=black] at (3.66, 4.66) {$1$};
     \node[draw,circle,inner sep=1.5pt,fill=orange, color=orange] at (3,6) {};
     
      \node [below, color=black] at (6, -0.2) {$e_1$};
       \node [below, color=black] at (0, -0.2) {$e_2$};

  \end{scope}
  
  
   \begin{scope}[shift={(9,-9.5)}]
       \draw [fill=yellow](0,0) -- (6,0)--(3,6)--cycle;
      \draw [thick, color=black](0,0) -- (6, 0);
     \draw [thick, color=black](0,0) -- (5.5, 1);
     \draw [thick, color=black](0,0) -- (5,2);
     \draw [thick, color=black](1.5, 3) -- (5,2);
     \draw [thick, color=black](1.5,  3) -- (4.33, 3.33);
     \draw [thick, color=black](1.5, 3) -- (3.66,  4.66);
  
     \draw [thick, color=black](0,0) -- (3,6);
     \draw [thick, color=black](6,0) -- (3,6);
    
      \node[draw,circle,inner sep=1.5pt,fill=black, color=black] at (6, 0) {};
     \node[draw,circle,inner sep=1.5pt,fill=black, color=black] at (0,0) {};    
               \node [right, color=black] at (0.2,0.5) {$4$};
     \node[draw,circle,inner sep=1.5pt,fill=black, color=cyan] at (5.5, 1) {};
                  \node [left, color=cyan] at (5.4, 1.3) {$7$};
     \node[draw,circle,inner sep=1.5pt,fill=red, color=black] at (5,2) {};
                 \node [left, color=black] at (4.5,2) {$3$};
     \node[draw,circle,inner sep=1.5pt,fill=blue, color=black] at (1.5, 3) {};
                 \node [right, color=black] at (1.5, 2.6) {$1$};
     \node[draw,circle,inner sep=1.5pt,fill=black, color=black] at (4.33, 3.33) {};
                  \node [left, color=black] at (4.3, 3.6) {$2$};
     \node[draw,circle,inner sep=1.5pt,fill=black, color=black] at (3.66, 4.66) {};
                  \node [left, color=black] at (3.66, 4.66) {$1$};
     \node[draw,circle,inner sep=1.5pt,fill=orange, color=orange] at (3,6) {};
     
      \node [below, color=black] at (6, -0.2) {$e_1$};
       \node [below, color=black] at (0, -0.2) {$e_2$};
       
         \draw [->, very thick, color=magenta](8 ,5)   --  (7, 5); 

  \end{scope}
  

   \begin{scope}[shift={(0, -9.5)}]
       \draw [fill=yellow](0,0) -- (6,0)--(3,6)--cycle;
      \draw [thick, color=black](0,0) -- (6, 0);
     \draw [thick, color=black](0,0) -- (5.5, 1);
     \draw [thick, color=black](0,0) -- (5,2);
     \draw [thick, color=black](1.5, 3) -- (5,2);
     \draw [thick, color=black](1.5,  3) -- (4.33, 3.33);
     \draw [thick, color=black](1.5, 3) -- (3.66,  4.66);
  
     \draw [thick, color=black](0,0) -- (3,6);
     \draw [thick, color=black](6,0) -- (3,6);
    
      \node[draw,circle,inner sep=1.5pt,fill=black, color=cyan] at (6, 0) {};
                  \node [left, color=cyan] at (6, 0.3) {$11$};
     \node[draw,circle,inner sep=1.5pt,fill=violet, color=violet] at (0,0) {};    
               \node [right, color=black] at (0.2,0.5) {$4$};
     \node[draw,circle,inner sep=1.5pt,fill=black, color=black] at (5.5, 1) {};
                  \node [left, color=black] at (5.4, 1.3) {$7$};
     \node[draw,circle,inner sep=1.5pt,fill=red, color=black] at (5,2) {};
                 \node [left, color=black] at (4.5,2) {$3$};
     \node[draw,circle,inner sep=1.5pt,fill=blue, color=black] at (1.5, 3) {};
                 \node [right, color=black] at (1.5, 2.6) {$1$};
     \node[draw,circle,inner sep=1.5pt,fill=black, color=black] at (4.33, 3.33) {};
                  \node [left, color=black] at (4.3, 3.6) {$2$};
     \node[draw,circle,inner sep=1.5pt,fill=black, color=black] at (3.66, 4.66) {};
                  \node [left, color=black] at (3.66, 4.66) {$1$};
     \node[draw,circle,inner sep=1.5pt,fill=orange, color=orange] at (3,6) {};
     
      \node [below, color=black] at (6, -0.2) {$e_1$};
       \node [below, color=black] at (0, -0.2) {$e_2$};
       
         \draw [->, very thick, color=magenta](8 ,5)   --  (7, 5); 

  \end{scope}


 \end{tikzpicture}
\end{center}
\caption{Computing the numerical decorations associated 
    to the vertices of the anthyphairetic  triangle $T(4, 11)$ starting from the 
    lotus $\Lambda(11/4)$.}  
   \label{fig:addvertweights}
    \end{figure}
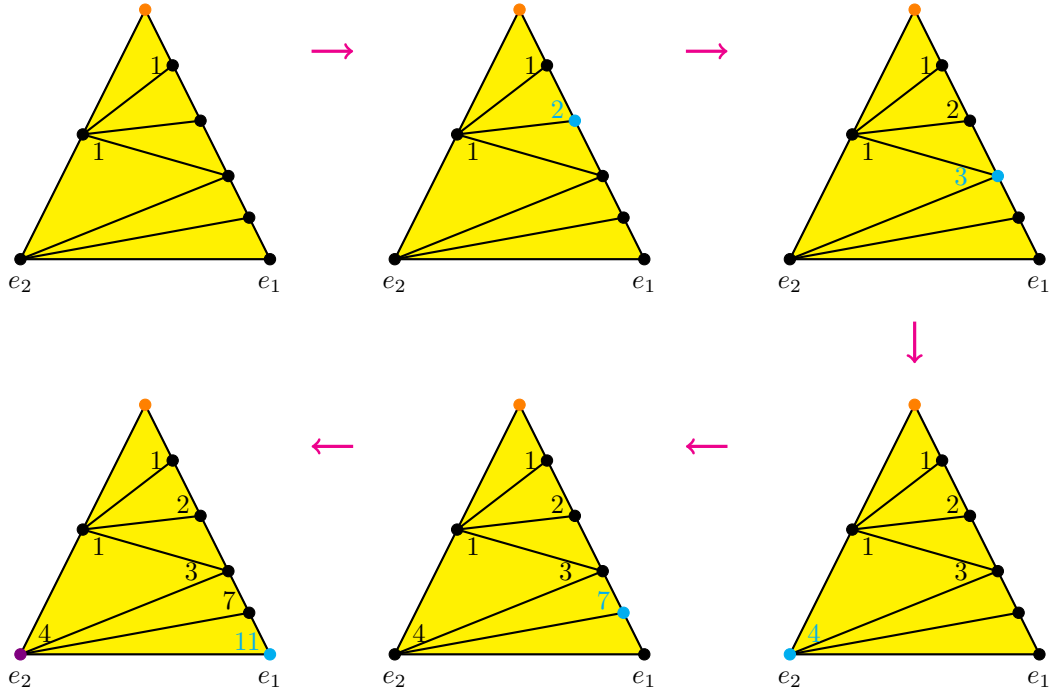

    \begin{proposition}  \label{prop:revalg}
     It is possible to pass graphically from the continued fraction expansion of a rational number 
  $b/a$ to the sequence of coprime pairs generated by the anthyphairetic process which starts from  
  $(a,b)$ by applying the following steps:   
 \begin{enumerate}
    \item[\bf S1)] Transform the continued fraction expansion of $b/a$ into the lotus $\Lambda(b/a)$. 
    \item[\bf S2)] Decorate both base vertices of the last petal of $\Lambda(b/a)$ by the integer $1$. 
    \item[\bf S3)] Descend the lotus petal by petal, attaching at each step the sum 
         of the two integers decorating the vertices of the petal to the third undecorated vertex.  
         This process is illustrated on Figure  \ref{fig:addvertweights}  for the pair $(4,11)$.
    \item[\bf S4)] Once this process is finished, construct the anthyphairetic rectangle $R(a,b)$ by 
      stacking the base edges of the petals of $\Lambda(b/a)$, endowed with the previous 
      numerical decorations, in the order in which they appear when one starts from the base of 
      $\Lambda(b/a)$.  For the pair $(4,11)$,  this step is the inverse of the right-hand arrow 
      of Figure \ref{fig:contrant}.
 \end{enumerate}
  \end{proposition}

 \begin{figure}[ht!]
     \begin{center}
\begin{tikzpicture}[scale=0.4]

   \begin{scope}[shift={(0, 0)}]
  
     \draw [fill=yellow](0,0) -- (2,3)--(4,-2)--cycle;
  
      \draw [thick, color=black](0,0) -- (2,3)--(4,-2)--cycle;
         
     \node[draw,circle,inner sep=1.5pt,fill=orange, color=black] at (0,0) {};
     \node[draw,circle,inner sep=1.5pt,fill=orange, color=black] at (2,3) {};
     \node[draw,circle,inner sep=2pt,fill=blue, color=cyan] at (4,-2) {};    
     
      \node [below] at (0,0) {$\beta$}; 
      \node [above] at (2,3) {$\alpha$}; 
      
        \draw [-> , very thick, color=magenta](5 , 2)   --  (6, 2);  
      
  \end{scope}


\begin{scope}[shift={(8, 0)}]
  
     \draw [fill=yellow](0,0) -- (2,3)--(4,-2)--cycle;
  
      \draw [thick, color=black](0,0) -- (2,3)--(4,-2)--cycle;
         
     \node[draw,circle,inner sep=1.5pt,fill=orange, color=black] at (0,0) {};
     \node[draw,circle,inner sep=1.5pt,fill=orange, color=black] at (2,3) {};
     \node[draw,circle,inner sep=2pt,fill=blue, color=cyan] at (4,-2) {};    
     
      \node [below] at (0,0) {$\beta$}; 
      \node [above] at (2,3) {$\alpha$}; 
      \node [below] at (4,-2) {\cyan{$\alpha + \beta$}}; 
      
  \end{scope}

   \end{tikzpicture}
\end{center}
  \caption{The basic operation to be performed inside the petals in order to apply the algorithm 
    of Proposition \ref{prop:revalg},  illustrated in Figure \ref{fig:addvertweights}.}  
  \label{fig:basicstep} 
    \end{figure}
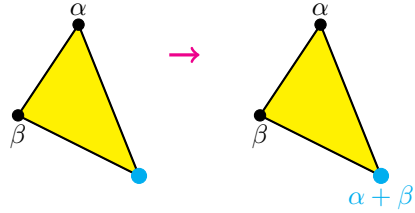     
    
     The proof that the algorithm of Proposition \ref{prop:revalg} 
 is correct is immediate once one notices that each petal gives rise to a subtraction rectangle, 
 in which the greatest of the two integers decorating the lower side is the sum of the two integers 
 decorating the upper side.

\begin{remark} 
    The previous algorithm is an example of the way one may use lotuses 
    as {\em computational architectures},  in the sense explained and amply illustrated 
    in the article \cite{GGP 26}. Namely, one may compute various invariants of 
    singularities by iterating very simple operations performed inside petals. Here, the 
    integers associated to the vertices of the lotus $\Lambda(b/a)$ (as illustrated on the 
    right side of Figure \ref{fig:contrant}) are the intersection numbers of the strict transform 
    of the binomial curve $C_{a,b}$ with the curves represented by those vertices, and the 
    operation to be iterated is shown in Figure \ref{fig:basicstep}. 
 \end{remark}


\bigskip 
\noindent
\textbf{\small{Authors' addresse:}}
\smallskip
\

\noindent
\small{P.\ Popescu-Pampu,
  Univ.~Lille, CNRS, UMR 8524 - Laboratoire Paul Painlev{\'e}, F-59000 Lille, France.
  \\
\noindent \emph{Email address:} \url{patrick.popescu-pampu@univ-lille.fr}}
\vspace{2ex}

\medskip
\end{document}